\documentclass[11pt]{amsart}
\usepackage{amsmath,amssymb,tabularx,graphicx,xypic}

\begin{document}

\bibliographystyle{amsalpha}

\title{Maximal Singular Loci of Schubert Varieties in $SL(n)/B$}

\date{\today}
\author{Sara C. Billey}\thanks{Work supported by NSF grant DMS-9983797}
\email{billey@math.mit.edu} 

\address{Author's address: 
Dept.\ of Mathematics, 2-363c\\ 
Massachusetts Institute of Technology\\ Cambridge, MA 02139}

\author{Gregory S. Warrington}
\email{gwar@math.harvard.edu} 

\address{Author's address: 
Dept.\ of Mathematics\\
Harvard University\\ Cambridge, MA 02138}



\renewcommand{\sectionmark}[1]{\markright{\thesection. #1}}

\theoremstyle{plain}
\newtheorem{thm}{Theorem}
\newtheorem{lem}[thm]{Lemma}
\newtheorem{cor}[thm]{Corollary}
\newtheorem{prop}[thm]{Proposition}
\newtheorem{conj}[thm]{Conjecture}

\theoremstyle{definition}
\newtheorem{dfn}[thm]{Definition}
\newtheorem{exmp}[thm]{Example}
\newtheorem{remark}[thm]{Remark}
\newtheorem{fact}[thm]{Fact}





\numberwithin{equation}{section}
\newcommand{\thmref}[1]{Theorem~\ref{#1}}
\newcommand{\corref}[1]{Corollary~\ref{#1}}
\newcommand{\propref}[1]{Proposition~\ref{#1}}
\newcommand{\lemref}[1]{Lemma~\ref{#1}}
\newcommand{\lemsref}[1]{Lemmas~\ref{#1}}
\newcommand{\factref}[1]{Fact~\ref{#1}}
\newcommand{\rmkref}[1]{Remark~\ref{#1}}
\newcommand{\defref}[1]{Definition~\ref{#1}}
\newcommand{\egref}[1]{Example~\ref{#1}}
\newcommand{\figref}[1]{Figure~\ref{#1}}
\newcommand{\figsref}[2]{Figures~\ref{#1},\ref{#2}}
\newcommand{\secref}[1]{Section~\ref{#1}}

\newcommand{\myfig}[2]{\begin{figure}[htbp]\begin{center}
      {\scalebox{.4}{\includegraphics{#1.eps}}}
      \caption{#2}\label{fig:#1}
    \end{center}\end{figure}}

\newcommand{\mysmfig}[2]{\begin{figure}[htbp]\begin{center}
      {\scalebox{.25}{\includegraphics{#1.eps}}}
      \caption{#2}\label{fig:#1}
    \end{center}\end{figure}}

\newcommand{\myfigtwo}[4]{\begin{figure}[htbp]\begin{center}\leavevmode
      {\scalebox{#1}{\includegraphics{#2.eps}}}
      {\scalebox{#1}{\includegraphics{#3.eps}}}
      \label{fig:#2}
      \caption{#4}
    \end{center}\end{figure}}

\newcommand{\myinline}[2]{{\scalebox{#1}{\includegraphics{#2.eps}}}}

\newcommand{\sfrac}[2]{\genfrac{\{}{\}}{0pt}{}{#1}{#2}}
\newcommand{\sumsb}[1]{\sum_{\substack{#1}}}  
        
\newcommand{\fsn}{\mathfrak{S}_n}
\newcommand{\fsp}[1]{S_{#1}}

\newcommand{\io}{{i_1}}
\newcommand{\iw}{{i_2}}
\newcommand{\inn}{{i_r}}

\newcommand{\bbR}{\mathbb{R}}
\newcommand{\bbC}{\mathbb{C}}
\newcommand{\bbZ}{\mathbb{Z}}
\newcommand{\bbN}{\mathbb{N}}
\newcommand{\bbQ}{\mathbb{Q}}

\newcommand{\wbi}{w^{\hat{i}}}
\newcommand{\xbi}{x^{\hat{i}}}
\newcommand{\zbi}{z^{\hat{i}}}
\newcommand{\wbio}{w^{\hat{1}}}
\newcommand{\wsbio}{ws^{\hat{1}}}
\newcommand{\xbio}{x^{\hat{1}}}
\newcommand{\zbio}{z^{\hat{1}}}

\newcommand{\tti}{t}
\newcommand{\wti}{\widetilde{w}}
\newcommand{\xti}{\widetilde{x}}
\newcommand{\uti}{\widetilde{u}}
\newcommand{\vti}{\widetilde{v}}
\newcommand{\yti}{\widetilde{y}}
\newcommand{\zti}{\widetilde{z}}

\newcommand{\defeq}{\overset{\text{def}}{=}}

\newcommand{\nts}{\negthickspace}

\newcommand{\trans}{\mathcal{T}}
\newcommand{\simref}{\mathcal{S}}


\newcommand{\msp}{\textsc{msp}}

\newcommand{\br}{\mathbf{a}}
\newcommand{\bal}{{\boldsymbol{\alpha}}}
\newcommand{\bbe}{{\boldsymbol{\beta}}}

\newcommand{\IH}{\operatorname{IH}}


\newcommand{\rp}[2]{\mathcal{R}({#1},{#2})}
\newcommand{\qp}[2]{\mathcal{E}_{{#1},{#2}}(x,w)}
\newcommand{\qpt}[2]{\mathcal{E}_{{#1},{#2}}(x,w)}

\newcommand{\rxw}{\mathcal{R}(x,w)}
\newcommand{\rytw}{\mathcal{R}(yt,w)}
\newcommand{\rysw}{\mathcal{R}(ys,w)}
\newcommand{\rysiw}{\mathcal{R}(ys_i,w)}
\newcommand{\qabxw}{\mathcal{E}_{a,b}(x,w)}
\newcommand{\qptxw}{\mathcal{E}_t(x,w)}
\newcommand{\qptabxw}{\mathcal{E}_{\tab}(x,w)}
\newcommand{\qbcxw}{\mathcal{E}_{b,c}(x,w)}
\newcommand{\qcdxw}{\mathcal{E}_{c,d}(x,w)}
\newcommand{\qbcyw}{\mathcal{E}_{b,c}(y,w)}

\newcommand{\qabxtwt}{\mathcal{E}_{a,b}(\widetilde{x},\widetilde{w})}
\newcommand{\qptxtwt}{\mathcal{E}_t(\widetilde{x},\widetilde{w})}
\newcommand{\qptabxtwt}{\mathcal{E}_{\tab}(\widetilde{x},\widetilde{w})}
\newcommand{\qbcxtwt}{\mathcal{E}_{b,c}(\widetilde{x},\widetilde{w})}
\newcommand{\qcdxtwt}{\mathcal{E}_{c,d}(\widetilde{x},\widetilde{w})}

\newcommand{\rxiwi}{\mathcal{R}(x^{-1},w^{-1})}
\newcommand{\rysiwi}{\mathcal{R}(sy^{-1},w^{-1})}
\newcommand{\ryw}{\mathcal{R}(y,w)}
\newcommand{\ryiwi}{\mathcal{R}(y^{-1},w^{-1})}
\newcommand{\rxtwt}{\mathcal{R}(\widetilde{x},\widetilde{w})}
\newcommand{\rxtw}{\mathcal{R}(xt,w)}

\newcommand{\digx}{\operatorname{mat}(x)}
\newcommand{\digxt}{\operatorname{mat}(\widetilde{x})}
\newcommand{\digwt}{\operatorname{mat}(\widetilde{w})}
\newcommand{\augdigx}{\operatorname{mat}'(x)}
\newcommand{\digw}{\operatorname{mat}(w)}
\newcommand{\digp}[1]{\operatorname{mat}(#1)}
\newcommand{\dpw}[1]{d_{#1,w}}
\newcommand{\duv}{d_{u,v}}
\newcommand{\dvw}{d_{v,w}}
\newcommand{\dzw}{d_{z,w}}
\newcommand{\dxw}{d_{x,w}}
\newcommand{\dyw}{d_{y,w}}
\newcommand{\txw}{\Theta_{x,w}}
\newcommand{\dxtwt}{d_{\widetilde{x},\widetilde{w}}}

\newcommand{\delxw}{\Delta(x,w)}
\newcommand{\delyw}{\Delta(y,w)}
\newcommand{\delytw}{\Delta(yt,w)}

\newcommand{\rxjw}{\mathcal{R}(x_i,w)}
\newcommand{\rxjmw}{\mathcal{R}(x_{i-1},w)}
\newcommand{\rxow}{\mathcal{R}(x_0,w)}

\newcommand{\area}{\mathcal{A}}
\newcommand{\pt}{\operatorname{pt}}

\newcommand{\ph}{\phi}
\newcommand{\pht}{\phi_t}
\newcommand{\phs}{\phi_s}
\newcommand{\phtp}{\phi_{t'}}
\newcommand{\phtal}{\phi_{t_{\alpha_j}}}

\newcommand{\fl}{\operatorname{fl}}
\newcommand{\ra}{\operatorname{unfl}}
\newcommand{\im}{\operatorname{Im}}
\newcommand{\slide}{\operatorname{\tau}}

\newcommand{\cP}{{\mathcal P}}
\newcommand{\xsing}{{X_w^{\text{sing}}}}
\newcommand{\maxsing}{\operatorname{maxsing}(X_w)}
\newcommand{\maxsingt}{\operatorname{maxsing}(X_{\widetilde{w}})}
\newcommand{\schub}[1]{X_{#1}}

\newcommand{\pxw}{P_{x,w}}
\newcommand{\psxw}{P_{sx,w}}
\newcommand{\pxsw}{P_{xs,w}}

\newcommand{\tpq}{t_{p,q}}
\newcommand{\tab}{t_{a,b}}
\newcommand{\talbt}{t_{\alpha,\beta}}
\newcommand{\tac}{t_{a,c}}
\newcommand{\tbc}{t_{b,c}}
\newcommand{\tcd}{t_{c,d}}
\newcommand{\tij}{t_{i,j}}
\newcommand{\tik}{t_{i,k}}
\newcommand{\tjk}{t_{j,k}}
\newcommand{\talj}{t_{\alpha_j}}
\newcommand{\tgd}{t_{\gamma,\delta}}

\newcommand{\ttpq}{t_{p,q}}
\newcommand{\ttab}{t_{a,b}}
\newcommand{\ttalbt}{t_{\alpha,\beta}}
\newcommand{\ttac}{t_{a,c}}
\newcommand{\ttbc}{t_{b,c}}
\newcommand{\ttcd}{t_{c,d}}
\newcommand{\ttij}{t_{i,j}}
\newcommand{\ttik}{t_{i,k}}
\newcommand{\ttjk}{t_{j,k}}
\newcommand{\ttalj}{t_{\alpha_j}}
\newcommand{\ttgd}{t_{\gamma,\delta}}

\newcommand{\sij}{s_{i_j}}
\newcommand{\sik}{s_{i_k}}


\begin{abstract}
  Schubert varieties in the flag manifold $SL(n)/B$ play a key role in
  our understanding of projective varieties.  One important problem is
  to determine the locus of singular points in a variety.  In 1990,
  Lakshmibai and Sandhya showed that the Schubert variety $X_w$ is
  nonsingular if and only if $w$ avoids the patterns $4231$ and
  $3412$.  In this paper we give an explicit combinatorial description
  of the irreducible components of the singular locus of the Schubert
  variety $X_w$ for any element $w\in\fsn$.  These irreducible
  components are indexed by permutations which differ from $w$ by a
  cycle depending naturally on a $4231$ or $3412$ pattern in $w$.  Our
  description of the irreducible components is computationally more
  efficient ($O(n^6)$) than the previously best known algorithms,
  which were all exponential in time.  Furthermore, we give simple
  formulas for calculating the Kazhdan-Lusztig polynomials at the
  maximum singular points.
\end{abstract}

\maketitle

\section{Introduction}
\label{sec:intro}

Schubert varieties play an essential role in the study of the
homogeneous spaces $G/B$ for any semisimple group $G$ and Borel
subgroup $B$; every closed subvariety in $G/B$ can be written as the
union of Schubert varieties, the classes of Schubert varieties form a
basis for the cohomology ring of $G/B$ and the Schubert varieties
correspond to the lower order ideals of a partial order associated to
$G/B$.  Specifically, this Bruhat order is an order on the $T$-fixed
points in $G/B$ where $T$ is the maximal torus in $B$.  The $T$-fixed
points, $e_w$, correspond bijectively with elements in the Weyl Group
$W = N(T)/T$ of $G$ and $T$.  A tremendous amount of information about
a Schubert variety can be obtained by examining the corresponding Weyl
group element.  Our main theorem gives a simple and efficient method
for giving the irrducible components of the singular locus of a
Schubert variety\footnote{While in the process of preparing this
submission, the authors have learned that Manivel \cite{manivel} has
recently independently proved a theorem equivalent to
\thmref{thm:main}.}.

In the late 1950's, Chevalley \cite{CHEV} showed that all Schubert
varieties in $G/B$ are nonsingular in codimension one.  Since that
time, many beautiful results on determining singular points of
Schubert varieties have surfaced (see \cite{BLak}).  By definition,
the Schubert variety $X_w$ is the closure of the $B$-orbit of $e_w$.
Therefore any point $p \in X_w$ is singular if and only if all points
in the orbit $Bp$ are singular.  Since the singular locus of a variety
is closed, the singular locus of $X_w$ is a union of Schubert
varieties indexed by the maximal elements $v<w$ such that $e_v$ is
singular in $X_w$.

Let $\maxsing$ denote the maximal set of Weyl group elements
corresponding to singular points in $X_w$ in Bruhat order, i.e.  $X_v$
is an irreducible component of the singular locus of $X_w$ if and only
if $v \in \maxsing$.  The goal of this paper is to give an explicit
algorithm for finding $\maxsing$ in the case where $G$ is
$SL_n(\mathbb{C})$, $B$ is the set of invertible upper triangular
matrices, $T$ is the set of invertible diagonal matrices, and $W$ is
the symmetric group $\fsn$.  The algorithm we present is very
efficient, $O(n^6)$, and removes the need to search through all
nonsingular $T$-fixed points (as is the case with previously known
techniques).

In type $A$ (i.e., $G=SL(n)$), smoothness is equivalent to rational
smoothness (\cite{Deodhar85}, see also \cite{carrkutt} in the case of
$ADE$) so the maximal singular locus of $X_w$ also determines the
maximal permutations $x \leq w$ for which the corresponding
Kazhdan-Lusztig polynomial is different from 1.  We use the explicit
form of $\maxsing$ to compute all Kazhdan-Lusztig polynomials at
maximal singular points (\msp's); they are either $1+q+\dots +q^{k}$ or
$1+q^{k}$ depending on whether the corresponding bad pattern is $4231$ or
$3412$ (respectively).

\section{Main results}\label{sec:main}

In 1990, Lakshmibai and Sandhya \cite{Lak-San} showed that the
Schubert variety $X_w \subset SL(n)/B $ is smooth at every point if
and only if the permutation matrix for $w$ does not contain any $4
\times 4$ submatrix equal to $3412$ or $4231$.  We use these two
permutation patterns to produce the maximal permutations below $w$
which correspond to points in the singular locus.  This verifies the
conjecture stated in \cite{Lak-San}.  (Gasharov, using a map similar
to the one we introduce in \secref{sec:phi}, shows in \cite{gash} that
the points constructed in \cite{Lak-San} are singular.  His result
proves one direction of this conjecture.)  In fact,
our proof starts from an arbitrary maximal singular $T$-fixed point
$e_x$ in $\schub{w}$ and shows that $w$ must contain a $4231$ or
$3412$ pattern and $x$ must contain a $2143$ or $1324$ pattern
(respectively).

The main theorem below shows that elements of $\maxsing$ are obtained
by acting on $w$ by certain cycles.  These cycles, described in the
following theorem, are best absorbed graphically in terms of the
permutation matrices $\digx$ and $\digw$.  Examples
are shown in Figures \ref{fig:mainthm1} and~\ref{fig:mainthm2}.

\begin{thm}\label{thm:main}  $X_x$ is an irreducible component of
  the singular locus of $X_w$ if and only if 
  \begin{equation*}
    \label{eq:4}
    x = w \circ (\alpha_1,\ldots,\alpha_m,\beta_k,\ldots,\beta_1)
  \end{equation*}
  for disjoint sequences
  \begin{align*}
    1&\leq \alpha_1 < \cdots < \alpha_m \leq n\text{, with } 
    w(\alpha_1) > \cdots > w(\alpha_m),\text{ and }\\
    1&\leq \beta_1 < \cdots < \beta_k \leq n\text{, with } 
    w(\beta_1) > \cdots > w(\beta_k),
  \end{align*}
  the interiors of the shaded regions in Figures \ref{fig:mainthm1}
  and \ref{fig:mainthm2} do contain any other 1's in the permutation
  matrix of $w$, and one of the following cases holds:
  \begin{enumerate}
  \item 4231 Case:
    \begin{gather*}
      k, m \geq 2\text{ and }\\
      \alpha_1 < \beta_1,\ldots, \beta_{k-1} < \alpha_2,\ldots, 
      \alpha_m < \beta_k \text{ and }\\
      w(\alpha_m) > w(\beta_1).
    \end{gather*}\label{item:10}
  \item 3412 Case:
    \begin{gather*}
      k,m \geq 2\text{ and }\\
      \beta_{1},\ldots , \beta_{k-1} < \alpha_1 < \beta_k < \alpha_2, \ldots, \alpha_k \text{ and }\\ 
      w(\alpha_{m-1}) > w(\beta_1) > w(\alpha_m) > w(\beta_2).
    \end{gather*}
    \label{item:1}
  \item 45312 Case:
    \begin{gather*}
      k = m = 2\text{ and }\\
      \beta_1 < \alpha_1 < \beta_2 < \alpha_2 \text{ and }\\        
      w(\alpha_1) > w(\beta_1) > w(\alpha_2) > w(\beta_2) \text{
      and }\\
    \text{ entries of $\digw$ in region A of \figref{fig:mainthm2}.2 are in
      decreasing order}.
    \end{gather*}
    \label{item:111}
  \end{enumerate}
\end{thm}
\myfig{mainthm1}{Example of Case \ref{item:10} of \thmref{thm:main}.
  $\circ$'s denote 1's in $\digw$, $\bullet$'s denote 1's in $\digx$.}
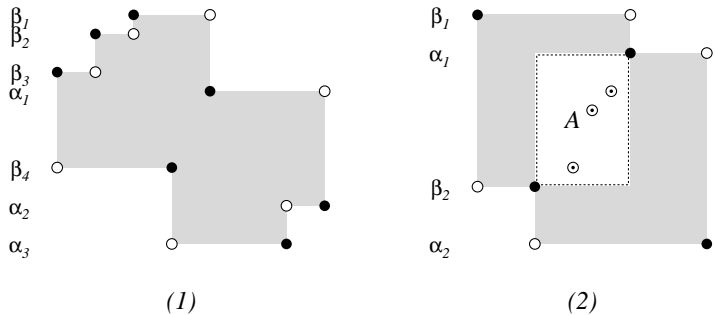
\myfig{mainthm2}{Examples for Cases \ref{item:1} and \ref{item:111} of
  \thmref{thm:main}.  For clarity in stating the theorem, region A in
  (2) is not shaded as it would be in the remainder of the
  paper. The $\odot$'s represent points where $\digx$ and $\digw$
  both have entries.}

After introducing basic notation in \secref{sec:prelim}, we then
introduce in \secref{sec:brf} the pictorial characterization of the
Bruhat order we rely on.  In sections \ref{sec:rxw} and \ref{sec:phi},
we discuss the Lakshmibai-Seshadri basis for the tangent space of a
Schubert variety indexed by transpositions and the set $\rxw = \{t: x <
xt \leq w\}$.  We also define a set of maps that allows us to relate
$\rxw$ and $\ryw$ when $x$ and $y$ differ by a transposition.  These
maps will then allow us to investigate not only whether a point $e_x$
is singular, but whether it is \emph{maximally} singular.  To describe
those permutations $x\in\maxsing$, we show that related permutations
$\xti$ must, among other qualities, avoid the patterns $231$, $312$
and $1234$.  We complete the description of $\maxsing$ in sections
\ref{sec:rest},\ref{sec:wtrest} and \ref{sec:maxsing}.  

The remaining sections contain applications arising from our
description of $\maxsing$.  In \secref{sec:lak}, we prove the
conjecture of Lakshmibai and Sandhya on the composition of $\maxsing$.
Using the tools we have developed, in \secref{sec:klp} we calculate
the values of the Kazhdan-Lusztig polynomials at maximal singular
points.  In \secref{sec:ex}, we give some example calculations
pertaining to the composition of $\maxsing$.  Finally, in
\secref{sec:ind}, we state a simple method for determining the number
of elements in $\maxsing $ in terms of pattern avoidance and
containment.

\section{Preliminaries}
\label{sec:prelim}

We begin by introducing our basic notation and terminology.  Let
$\fsn$ denote the symmetric group on $n$ letters. 
We will view elements of $\fsn$ as permutations on
$[1,\ldots,n]$.  To this end, we identify $s_i$ with the
transposition $(i,i+1)$.  Let $w(i)$ be the image of $i$ under the
permutation $w$.  We have a one-line notation for a permutation
$w$ given by writing the image of $[1,\ldots,n]$ under the action of
$w$: $[w(1),w(2),\ldots,w(n)]$.  We will also often utilize the
permutation matrix for $w$ (denoted $\digw$).

We use the standard presentation
\begin{equation}
  \begin{split}
    \fsn = \langle s_1,\ldots,s_{n-1} :\ &s_i^2 = 1,\\
    &s_is_j = s_j s_i \text{ for } |i-j| > 1 \text{, and }\\
    &s_is_{i+1} s_i = s_{i+ 1} s_i s_{i+1}\rangle.
  \end{split}
\end{equation}
Let $\simref =\{s_i\}_{i\in [1,\ldots,n-1]}$ denote the generating set for
$\fsn$.  Let $\trans$ denote the set of all \textit{transpositions}
in $\fsn$.  The elements of $\trans$ are all the conjugates of elements in
$\simref$:
\begin{equation}
 \trans = \{t_{j,k} = s_js_{j+1}\cdots s_{k-2}s_{k-1}s_{k-2}
     \cdots s_{j+1} s_j: 1 \leq j < k \leq n\}.   
\end{equation}
If we wish to refer to a transposition $t$ that affects positions $a$ and
$b$, but the relative values of $a$ and $b$ are unknown, we will
write $t_{\{a,b\}}$.  

An \textit{expression} is any product of generators $s_i$.  The
\textit{length} $l(w)$ of an element $w\in \fsn$ is the
minimum $r$ for which we have an expression $w = s_\io\cdots s_{i_r}$.
A \textit{reduced expression} $w = s_{i_1}\cdots s_{i_r}$ is an
expression for which $l(w) = r$.
It is a standard fact that
\begin{equation}
  l(w) = \#\{1 \leq i < j \leq n: w(i) > w(j)\}.
\end{equation}

\begin{dfn} Let $x,w \in \fsn$, $p,q\in \bbZ$.
  Define the \textit{rank function} for $w$ by $r_w(p,q)= \#\{i
  \leq p: w(i) \geq q\}$.  Also, the \textit{difference function} for the
  pair $x,w$ is defined by $\dxw(p,q) = r_w(p,q) - r_x(p,q)$.
\end{dfn}

In $SL(n)/B$, Schubert varieties can be defined in terms of the rank
function.  To do this, fix a basis $\langle u_1,u_2,\cdots,u_n \rangle$ of
$\bbC^n$.  This fixes a base flag $(F_1 \subset F_2 \subset \cdots
\subset F_n)$ where $F_i$ is the span of $\langle
u_n,\ldots,u_{n-i+1}\rangle$.   We identify $SL(n)/B$ with the
set of all complete flags of vector spaces $(V_1 \subset V_2 \subset
\cdots \subset V_n \cong \bbC^n)$ with $\dim V_i = i$.  Then
\begin{equation}\label{e:schubdef} X_w = \{( V_1 \subset V_2 \subset
\cdots \subset V_n): \dim(V_p \cap F_{n-q+1}) \leq r_w(p,q)\}.
\end{equation} 
(This definition is equivalent to that given in \cite{Fulton-book}.)
Note that the flag
\begin{equation}
  e_w = (\langle u_{w(n)} \rangle \subset \langle u_{w(n)},u_{w(n-1)}
  \rangle \subset \cdots \subset \langle u_{w(n)},\ldots,u_{w(1)} \rangle)  
\end{equation}
is an element of $X_w$.  Furthermore, $e_w$ is fixed by the left $T$ action.
Also, since $X_{v}$ can be viewed as the closure of the orbit
$Be_{v}$, $X_{v} \subseteq X_{w}$ if and only if $e_{v} \in X_{w}$.
This defines a partial order, called the Bruhat (or Bruhat-Chevalley)
order on $\fsn $ by
\begin{equation}
  v \leq w \Leftrightarrow X_v \subseteq X_w.
\end{equation}

The Bruhat order has a number of characterizations (see, e.g.,
\cite{Hum}).  One of the most common definitions is as the transitive
closure of the relations $vt < v$ for $t\in\trans$ if $l(vt) < l(t)$.
However, we prefer to work with a more graphical characterization
which follows directly from the definition of the rank and difference
functions above.  The corresponding ``Bruhat pictures'' that we
associate to each pair $x \leq w$ will be discussed in the next
section.  These pictures will rely on the below two conclusions.

 \begin{lem}\label{lem:bdf} We have $x \leq w$ if and only if $\dxw$
is everywhere non-negative.  \end{lem}

\begin{cor}\label{lem:trans}
  If $x \leq y \leq w$, then $\dxw - \dyw$ is everywhere non-negative.
\end{cor}

The Bruhat graph of $w$ is the graph with vertices labeled by $\{v
\leq w\}$ and $v_1$ is joined to $v_2$ by a directed edge if $v_1=v_2
t$ for some $t \in \trans$ and $v_1<v_2 $ in Bruhat order.  This graph
plays a central role in the study of Schubert varieties.  For example,
Lakshmibai and Seshadri have shown that in $SL(n)/B$, the tangent
space to $X_w$ at $e_x$ has a basis indexed by $\{t \in \trans : xt
\leq w\}$, i.e. the edges of the Bruhat graph adjacent to $x$.  This
fact forms the main criterion we will use in Section~\ref{sec:rxw} for
smoothness at a point.  In fact, since $xt < x$ implies $xt < w$ we
will just need to consider the edges ``going up'' from $x$ in the
Bruhat graph of $w$.  As mentioned in the introduction, this set will
be denoted by $$\rxw := \{t\in \trans: x < xt \leq w\}.$$

Over the last few years, it has become apparent that properties of the
Bruhat order can often be efficiently characterized by ``pattern
avoidance'' \cite{BL,B,hexagon,LS2,Stem4}.  We say that
$w = [w(1),\ldots,w(n)]$ \textit{avoids the pattern} $v =
[v(1),\ldots,v(k)]$ for $k \leq n$ if we cannot find $1 \leq i_1 <
\cdots < i_k \leq n$ with $w(i_1),\ldots,w(i_k)$ in the same relative
order as $v(1),\ldots,v(k)$ - i.e., no submatrix of $\digw$ on rows
$i_1,\ldots,i_k$ and columns $w(i_1),\ldots,w(i_k)$ is the permutation
matrix of $v$.  Our characterization of the maximal singular locus is
yet another example of the utility of this notion.

More generally, we can define pattern avoidance or containment in
terms of the following flattening function.  For any set $Z = \{z_1 <
z_2 < \cdots < z_k\} \subseteq [1,\ldots,n]$, and $x\in \fsn$, define
$\fl_Z(x)$ to be the ``flattened'' permutation on $[1,\ldots,k]$ whose
elements are in the same relative order as $[x(z_1), \ldots, x(z_k)]$.
When the set $Z$ is clear from context, we will abbreviate $\fl_Z(x)$
by $\overline{x}$.  We will also write $\fl(i,j,\ldots,k)$ for the
flattened permutation on the sequence $i,j,\ldots, k$ and write $\xbi$
for $\fl_{[1,\ldots,n] \setminus \{i\}}(x)$.

It will also be useful to have notation for an ``unflattening''
operator.  Given a permutation $x\in\fsn$, a set $Z\subseteq
[1,\ldots,n]$, and a permutation $u\in\fsp{k}$, we can define a new
permutation $\ra^x_{Z}(u)\in\fsn$ by requiring that
\begin{enumerate}
\item $\fl_Z(\ra^x_{Z}(u)) = u$, and
\item $x(a) = (\ra^x_{Z}(u))(a)$ if $a\in[1,\ldots,n] \setminus Z$.
\end{enumerate}
When $x$ and $Z$ are clear from context, we abbreviate $\ra^x_{Z}(u)$
by $\widehat{u}$.
\begin{exmp}
  For $x = [5,2,4,1,6,3]$ and $Z = \{3,5,6\}$, we have $\fl_Z(x) =
  [2,3,1]$ and $(\ra^x_{Z})([3,1,2]) = [5,2,6,1,3,4]$.  Note that $x =
  \widehat{\fl_Z(x)}$.
\end{exmp}

\section{Bruhat pictures}
\label{sec:brf}

Our main theorem is concerned not only with determining which points
in a given Schubert variety are singular, but which are
\emph{maximally} singular.  The function $\dxw$ affords us a graphical
view of the Bruhat order.  Most importantly, it lets us see the
set $\rxw$.  We will now introduce the graphical notation utilized in
the remainder of the paper that allows us to do this.  A diagram
displaying the notation we are about to describe is offered in
\figref{fig:diagex}.  

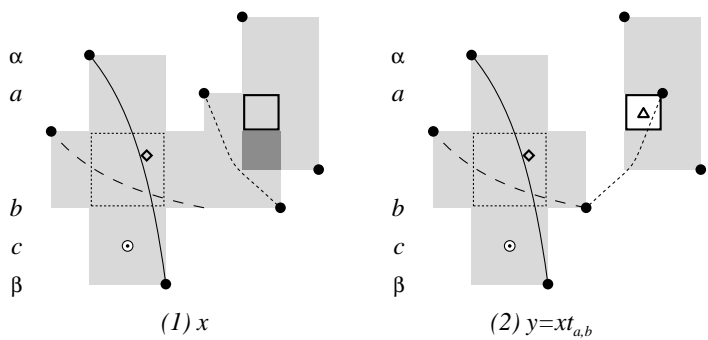
\myfig{diagex}{We see (among other facts) that $\dxw(\diamond) \geq
  1$, $\dyw(\triangle) = \dxw(\triangle) - 1 \geq 0$, $\pt_x(c) =
  \pt_w(c)$ and $\talbt\in\rxw$.}

First, we plot, as black disks, all or some of the positions containing
1's in the permutation matrix $\digx$ of $x$.  We will sometimes
overlay $\digx$ and $\digw$.  In these cases, 1's in $\digw$ will
be marked by open circles.  Points that are simultaneously in both
diagrams will consist of a black disk and a larger concentric circle.
Let $[a,b]\times [c,d]$ denote the set of all points $(p,q)\in\bbR^2$
such that $a \leq p \leq b$ and $c\leq q\leq d$.  The following
notation will be handy:
\begin{dfn}
  For $\tpq\in\rxw$, set 
  \begin{align}
    \area_{p,q} &:= \area_{p,q}(x) = [p+1,q-1]\times [x(p)+1,x(q)-1],\\
 \overline{\area_{p,q}} &:= \overline{\area_{p,q}(x)} = [p,q]\times [x(p), x(q)],\\
   \pt_x(c) &:= (c,x(c)) \text{ for } c\in[1,\ldots,n].
  \end{align}
\end{dfn}
Along with the points of $\digx$, we will often shade parts of our
diagram in order to specify that $\dxw$ satisfies a particular
inequality on a given region.  Light shading on a region signifies
that $\dxw \geq 1$ on that region.  Dark shading signifies $\dxw \geq
2$.  No shading places no restrictions on the values $\dxw$.  A region
with a black border is one where $\dxw$ achieves the minimum possible
value allowed by the shading on that region.  Dotted borders are used
to demarcate regions we wish to discuss in the text.

As mentioned above, the great utility of these diagrams arises from
being able to visualize $\rxw$ along with the information on the
Bruhat order.  To see how we do this, suppose we have some reflection
$\tab\in\rxw$ (which implies $x < x\tab\leq w$).  Now compare the shading
(with respect to $w$) in $\digp{x}$ and $\digp{x\tab}$.  We see (as in
\figref{fig:diagex}), that in the region $\area_{a,b}(x)$,
$\dpw{x\tab} = \dxw - 1$.  Hence, by \lemref{lem:bdf}, we can state
the following:
\begin{fact}\label{fact:shady}
  Let $\tab\in \trans$ with $x < x\tab$.  The transposition $\tab$ is in
  $\rxw$ if and only if it corresponds to a region in $\digx$ that is
  entirely shaded (i.e., $\dxw|_{\area_{a,b}} \geq 1$).  An example is
  given in \figref{fig:shadedex}
\end{fact}
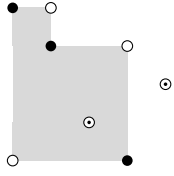
\myfig{shadedex}{In the above configuration, $\rxw = \{t_{1,2}, t_{2,4},
  t_{2,5}, t_{4,5}\}$.}  

Note that the values of $\dxw$ on the region $\overline{\area_{a,b}}\ 
\setminus \area_{a,b}$ are not considered in determining the
membership of $\tab$ in $\rxw$.
  
In order to highlight reflections that we are particularly interested
in, we will often draw an arc in our diagram.  A solid or dotted curve
connecting two points in $\digx$ will denote an element of $\rxw$.  A
dotted curve will be used to designate $t$ when we are particularly
interested in $y = xt$.  A dashed curve will be used when we wish to
mark a reflection $t'\in\ryw$.  Of course, if $tt' \neq t't$, and our
picture is of $\digx$, then only one of the endpoints of our dashed
curve will correspond to a point of $\digx$.  There are numerous
instances when an arc corresponds to an element of $\rxw \cap \ryw$.  In
this case, whether we use a solid or dashed arc depends on context.

The following lemma will be used several times in future sections.  It
allows us to infer the presence of points in $\digx$ in a region based
on a particular common pattern of shading.

\begin{lem}\label{lem:rhombus} 
  Let $x < w$ and suppose $p,p',q,q'\in \bbZ$ such that 
  \begin{enumerate}
  \item $p < p'$, $q < q'$,
  \item $\dxw(p,q') = 0$,
  \item $\dxw(p,q) = \alpha$, $\dxw(p',q') = \beta$,
    $\dxw(p',q) = \gamma$.
  \end{enumerate}
  Then there exist at least $\alpha+\beta-\gamma$ values $m$ such that
  $\pt_x(m) \in [p+1,p'-1] \times [q+1,q'-1]$ with $x(m) \neq w(m)$.
\end{lem}

\begin{proof}
  Let's define four regions as follows:
  \begin{itemize}
  \item $A = [1,p]\times [q,q'-1]$,
  \item $B = [p+1,p']\times [q',n]$,
  \item $C = [p+1,p']\times [q,q'-1]$,
  \item $D = [1,p]\times [q',n]$.
  \end{itemize}

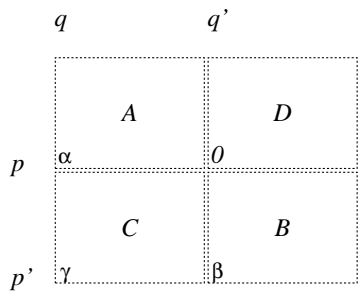
  \myfig{rhombus}{We have indicated certain values of $\dxw$
    at the lower left corner of each region.}

  For every subset $R \subset [1,n]\times [1,n]$, define
  \begin{equation}
    \label{deldfn}
    \txw(R) = \#\{(p,q)\in R: q = w(p)\} -
    \#\{(p,q)\in R: q = x(p)\}.
  \end{equation}
  Then
  \begin{align*}
    \dxw(p,q') &= 0 \text{ implies that } \txw(D) = 0,\\    
    \dxw(p',q') &= \beta \text{ implies that } \txw(B) = \beta,
    \text{ and }\\    
    \dxw(p,q') &= \alpha \text{ implies that } \txw(A) = \alpha.    
  \end{align*}
  Now, 
  \begin{align}
    \dxw(p',q) &= \txw(A) + \txw(B) + \txw(C) + \txw(D), \text{ so }\\
                \gamma &= \alpha + \beta + \txw(C) + 0.
  \end{align}
  So $\txw(C) = -(\alpha + \beta - \gamma)$ and there are
  exactly $\alpha + \beta - \gamma$ more 1's of $\digx$ than 1's of
  $\digw$ in region $C$.  This finishes the proof.
\end{proof}

\section{A criterion for maximal smoothness}
\label{sec:rxw}

To prove Theorem~\ref{thm:main}, we start from the fact that (by
definition) $X_w$ is smooth at $e_x$ if and only if the dimension of
the Zariski tangent space at that point is equal to $l(w) = \dim(X_w)$.
Lakshmibai and Seshadri, \cite{LakSes}, describe the dimension of this
tangent space in terms of the root system.  Using the fact that
$\#\{t\in \trans: xt < x\} = l(x)$, we can paraphrase their result as:

\begin{thm}\cite{LakSes}\label{fact:basic}
  The Schubert variety $X_w \in SL (n)/B$ is smooth at $e_x$ if and
only if $\#\rxw := \#\{t\in \trans: x < xt \leq w\}$ equals $l(w) - l(x)$.
\end{thm}
This yields the following characterization of the permutations
in $\maxsing$:
\begin{fact}\label{fact:singcrit}
  $x\in\maxsing$ if and only if
  \begin{enumerate}
  \item $\#\rxw > l(w) - l(x)$ and
  \item for all $t\in\rxw$, $\#\rp{xt}{w} = l(w) - l(xt)$.
  \end{enumerate}
\end{fact}

As may be ascertained from \thmref{thm:main}, the criteria for $x$ to
be an element of $\maxsing$ are local in nature.  This implies that we
may concentrate on only certain indices in our permutation $w$ in
order to determine $\maxsing$.  We now describe these indices
explicitly.

\begin{dfn}\label{dfn:tilde}
  Let 
  \begin{equation}
    \delxw = \{i, 1 \leq i \leq n: \exists\ j, 1\leq j \leq n,
    \text{ with } t_{\{i,j\}}\in\rxw\}.
  \end{equation}
  For $\delxw = \{d_1 < d_2 < \cdots < d_k\}$,
  set 
  \begin{align}
    \xti &= \fl([x(d_1),x(d_2),\ldots,x(d_k)]) \text{ and }\\
    \wti &= \fl([w(d_1),w(d_2),\ldots,w(d_k)]).
\end{align}
\end{dfn}
Note that $\xti$ and $\wti$ are permutations in $\fsp{k}$.

We state here for reference the following useful characterization of $\delxw$.
\begin{cor}\label{cor:deliff}
  Let $x\leq w$.  Then $\dxw(\pt_x(p)+(-1,0)) = 0$
  and $\dxw(\pt_x(p)+(0,1)) = 0$ if and only if
  $p\not\in\delxw$.
\end{cor}

We now give a sufficient condition for an index $b$ to be in $\delxw$.
\begin{prop}\label{prop:lblf}
  Suppose $x < w$ and $x(b) \neq w(b)$ with $1 \leq b \leq n$.
  \begin{enumerate}
  \item If $w(b) < x(b)$, then $\exists\ a < b$ with $\tab\in\rxw$
    and $x(a) \neq w(a)$.
  \item If $w(b) > x(b)$, then $\exists\ c > b$ with $\tbc\in\rxw$
    and $x(c) \neq w(c)$.
  \end{enumerate}
\end{prop}
\myfig{lblf}{$\diamond = (p,q')$ and $\triangle = (p',q)$.}

\begin{proof}
  First we prove the case of $w(b) < x(b)$.  Note that 
  \begin{equation}
    \dxw(b-1,x(b)) = 1 + \dxw(b,x(b)) \geq 1   
  \end{equation}
  since $w(b) < x(b)$.  Let $p' = b-1$.  Choose $q$ as large as
  possible such that $q < x(b)$ and $\dxw(p',q) = 0$ (see
  \figref{fig:lblf}).  Such a $q$ must exist since $\dxw(\cdot,0) =
  0$.  Now choose $p$ as small as possible such that $p < p'$ and
  $\dxw(g,h) \geq 1$ for all $g,h$ with
  $(g,h)\in[p+1,p']\times[q+1,x(b)]$.  Then there exists a $q'$, $q <
  q' \leq x(b)$ such that $\dxw(p,q') = 0$.  By construction,
  \begin{equation*}
    \dxw(p,q') = 0,\ \dxw(p',q) = 0,\ \dxw(p,q) \geq 0 \text{ and }\dxw(p',q')
    \geq 1.
  \end{equation*}
  That is (in the notation of \lemref{lem:rhombus}), $\alpha \geq 0$,
  $\beta \geq 1$ and $\gamma = 0$.  So by this lemma, there exists an
  $a$ such that $\pt_x(a)\in [p+1,p']\times [q,q'-1]$ and $x(a) \neq
  w(a)$.  Then $\dxw|_{\area_{a,b}} \geq 1$, so by
  \factref{fact:shady}, $\tab\in\rxw$.  This proves our claim.

  To prove the case of $w(b) > x(b)$, it is easiest to use dual
  rank and difference functions: 
  \begin{align}
    r'_w(p,q) &:= \#\{i \geq p: w(i) \leq q\},\\
    d'_{x,w} &:= r'_w - r'_x.\label{eq:dualdxw}
  \end{align}
  One can check that $x \leq w$ if and only if $d'_{x,w} \geq 0$ and
  then argue as above using this new rank function.  (Note, to define
  $X_w$ using $r'_w$, we have to modify our fixed flag $F_\cdot$ and 
  (\ref{e:schubdef}).)
\end{proof}

\begin{cor}\label{cor:uprfln}
  If $x \leq w$ and $\dxw(\pt_x(b)) > 0$, then there exists $b' < b$
  with $t_{b',b}\in\rxw$.
\end{cor}

\propref{prop:lblf} tells us that if $x(i) \neq w(i)$ then
$i\in\delxw$.  It turns out that the question of whether or not
$x\in\maxsing$ depends only on the pair $\xti,\ \wti$.  This is borne
out by the following simple facts.  They will be used without comment
in the remainder of the paper.

\begin{lem}\label{lem:comb}
  We have the following:
  \begin{enumerate}
  \item If $x(i) = w(i)$, then $\xbi \leq \wbi \Longleftrightarrow
    x \leq w$.  
  \item $\xti \leq \wti \Longleftrightarrow x \leq w$.
  \end{enumerate}
\end{lem}
\begin{proof}
  The first equivalence follows from \lemref{lem:bdf} by comparing
  $\dxw$ and $d_{\xbi,\wbi}$.  The second follows from the first by
  noting that $i\in\delxw$ whenever $x(i) \neq w(i)$.
\end{proof}

\begin{prop}\label{prop:rxtwt}
  We have the following:
  \begin{enumerate}
  \item $l(w) - l(x) = l(\wti) - l(\xti)$.
  \item There exists a bijection $\rxtwt \overset{\sim}{\longrightarrow} \rxw$.
  \item $x\in\maxsing$ if and only if $\xti \in\maxsingt$.
  \end{enumerate}
\end{prop}
\begin{proof}
  Pick some $i\not\in\delxw$.  Now,
  \begin{align}
    l(w) - l(x) - (l(\wbi) - l(\xbi)) = \dxw(\pt_x(i)) + \dxw'(\pt_x(i)).
  \end{align}
  We know by \corref{cor:uprfln} that $\dxw(\pt_x(i)) = 0$.  Applying
  \lemref{lem:rhombus} with $p = i,\ p' = n,\ q = 1$ and $q' = x(i)$,
  we see that $\dxw'(\pt_x(i))$ is also 0.  This proves Part 1.  Part
  2 follows immediately from \factref{fact:shady} and the definition
  of $\delxw$ by comparing $\dxw$ and $\dxtwt$.  Part 3 follows from
  the first two parts along with \corref{cor:delxy} (stated below). 
\end{proof}

\section{The Map $\pht$}
\label{sec:phi}

In \factref{fact:singcrit} we claimed that $\maxsing$ can be
identified in terms of $\rxw$ for $x \leq w$.  To do this, we will
need to relate $\rxw$ to $\ryw$ when $x,y$ differ by an element of
$\trans$.  So, for every triple $yt < y \leq w$ with $t\in \trans$, we
will define a map $\ph_{t}^{y,w}: \ryw \longrightarrow \trans$.  In
\thmref{thm:phiinj} we will show that the image is actually contained
in $\rytw$.  The values of $y,w$ are usually clear from context and we
will often abbreviate $\ph_{t}^{y,w}$ as $\pht$.  Propositions
\ref{prop:phimax},~\ref{prop:phimaxl} and ~\ref{prop:reci} contain
interesting implications of the map $\pht$.

A similar map has been defined by Gasharov \cite{gash} for the purpose
of showing that certain elements constructed by Lakshmibai and Sandhya
in \cite{Lak-San} are, in fact, singular points.  See \secref{sec:lak}
for details.

\begin{dfn}\label{def:phi} 
  Fix $yt < y \leq w$.  Given some $t'\in\ryw$, if $t$ and $t'$
  commute, we define $\pht(t') = t'$.  Otherwise, we can find $a < b <
  c$ such that $d\not\in \{a,b,c\}$ implies $y(d) = yt(d) = yt'(d)$.
  Then we define $\pht^{y,w}(t')$ according to Table ~\ref{tab:dfn}.
  \begin{table}
    \setlength{\extrarowheight}{4pt}
  \begin{tabular}[c]{|l|>{$}c<{$}>{$}c<{$}>{$}c<{$}>{$}c<{$}!{\vrule
        width 2pt}>{$}c<{$}!{\vrule width 2pt}
      >{$}c<{$}>{$}c<{$}|}\hline
    Case & \overline{y} & t & t' & \tac\in\ryw &
    \pht(t') & \overline{yt'} & \overline{yt\pht^{y,w}(t')}\\\hline
    $A.i)$ & 213 & \tab & t_{a,c} & \checkmark & t_{b,c} & 312 & 132\\
    $\ \ \ ii)$ & 213 & \tab & t_{b,c} & \times & t_{b,c} & 231 & 132\\
    $\ \ \ iii)$ & 213 & \tab & t_{b,c} & \checkmark & t_{a,c} & 231 & 321\\
    \hline

    $B.i)$ & 132 & \tbc & t_{a,c} & \checkmark & t_{a,b} & 231 & 213\\
    $\ \ \ ii)$ & 132 & \tbc & t_{a,b} & \times & t_{a,b} & 312 & 213\\
    $\ \ \ iii)$ & 132 & \tbc & t_{a,b} & \checkmark & t_{a,c} & 312 & 321\\
    \hline

    $C.i)$ & 312 & \tab & t_{b,c} & \times & t_{a,c} & 321 & 231\\
    \ \ \ $ii)$ & 312 & \tac & t_{b,c} & \times & t_{b,c} & 321 & 231\\\hline

    $D.i)$ & 231 & \tbc & t_{a,b} & \times & t_{a,c} & 321 & 312\\
    \ \ \ $ii)$ & 231 & \tac & t_{a,b} & \times & t_{a,b} & 321 & 312\\\hline
  \end{tabular}

  \caption{Definition of map $\pht^{y,w}$.  We have split into cases
    indexed by $\fl_{abc}(y)$ and whether $t' = t_{a,b}$,
  $t' = t_{a,c}$ or $t' = t_{b,c}$.  Note that the matter of 
    inclusion of
  $\tac$ in $\ryw$ is determined by the first three columns in
  Cases A.i,B.i,C and D.  The final two columns are used in proving
  that $\pht$ maps $\ryw$ into $\rytw$.}
  \label{tab:dfn}
\end{table}
\end{dfn}

\myfig{phidfn1}{Graphical depiction of Case A from \defref{def:phi}.
  The dashed (resp. dotted, solid) arcs represent $t'$ (resp. $t$,
  $\pht(t')$).}

\myfig{phidfn2}{Graphical depiction of Case C from \defref{def:phi}.
  The dashed (resp. dotted, solid) arcs represent $t'$ (resp. $t$,
  $\pht(t')$).}

\begin{remark}
  It is not sufficient to define $\pht(\tbc) = t_{\{t(b),t(c)\}}$.
  For example, in the situation of Case C.ii, where $t = \tac$ and $t'
  = \tbc$, we have $t_{\{t(b),t(c)\}} = t_{a,b}\not\in\rp{yt}{w}$.
\end{remark}

\begin{exmp}
  In order to elucidate the definition, we give here several example of the
  map $\ph_t$.  Let $w = [2,4,5,3,1]$ and $x_0 = [2,1,5,4,3]$.  We see
  that $\rp{x_0}{w} = \{t_{2,4},t_{2,5}\}$ (see \figref{fig:phiex}.1).
  \myfig{phiex}{}
  Let $x_1 = x_0t_{1,2} < x_0$.  We see that
  \begin{enumerate}
  \item $\ph_{t_{1,2}}(t_{2,4}) = t_{2,4}$ (Case A.ii).
  \item $\ph_{t_{1,2}}(t_{2,5}) = t_{2,5}$ (Case A.ii). 
  \end{enumerate}
  Now, $\rp{x_1}{w} = \{t_{1,2},t_{2,4},t_{2,5}\}$ (see 
  \figref{fig:phiex}.2).  Let $x_2 = [1,2,5,3,4] = x_1t_{4,5} < x_1$.  We see that
  \begin{enumerate}
  \item  $\ph_{t_{4,5}}(t_{1,2}) = t_{1,2}$ (Because $t_{4,5}t_{1,2} = t_{1,2}t_{4,5}$).
  \item  $\ph_{t_{4,5}}(t_{2,4}) = t_{2,5}$ (Case B.iii).
  \item  $\ph_{t_{4,5}}(t_{2,5}) = t_{2,4}$ (Case B.i).
  \end{enumerate}
  Now $\rp{x_2}{w} = \{t_{1,2}, t_{2,4}, t_{2,5}, t_{4,5}\}$.  Let
  $x_3 = [1,2,3,5,4] = x_2t_{3,4} < x_2$ (see \figref{fig:phiex}.3).  We see that
  \begin{enumerate}
  \item  $\ph_{t_{3,4}}(t_{1,2}) = t_{1,2}$ (Because $t_{3,4}t_{1,2} = t_{1,2}t_{3,4}$).
  \item  $\ph_{t_{3,4}}(t_{2,4}) = t_{2,3}$ (Case B.i).
  \item  $\ph_{t_{3,4}}(t_{2,5}) = t_{2,5}$ (Because $t_{3,4}t_{2,5} = t_{2,5}t_{3,4}$).
  \item  $\ph_{t_{3,4}}(t_{4,5}) = t_{3,5}$ (Case C.i).
  \end{enumerate}
\end{exmp}

\begin{remark}
  It is possible for $\#\rytw > \#\ryw + 1$ for $yt < y \leq w$ and
  $t\in \trans$.  For example, let $w = [4,2,3,1]$, $y =
  [2,4,1,3]$ and $x = yt_{2,3} = [2,1,4,3]$.  Then $\rp{y}{w} =
  \{t_{1,2}, t_{3,4}\}$ and $\rytw = \{t_{1,3},t_{1,4},t_{2,3},t_{2,4}\}$.
\end{remark}

It is clear from the definition that $\pht(t')$ is always a reflection
and that $yt < yt\cdot \pht(t')$.  But for $\pht$ to be useful, we
will need the following property.


\begin{thm}\label{thm:phiinj} 
  Fix $yt < y \leq w$.  The map $\pht(\ryw) \hookrightarrow \rytw
  \setminus \{t\}$ is injective.
\end{thm}

For the proof of this theorem, we'll need the following simple lemma.

\begin{lem}\label{lem:z}
  Let $u,w\in\fsn$ and suppose $1\leq i < j < k \leq n$ such that
  $\fl_{ijk}(u) = 123$.  If both
  \begin{align}
    w &\geq x = u\circ (k,j,i) \text{ (i.e., } \overline{x} = 312\text{) and }
    \label{eq:1}\\
    w &\geq y = u\circ (i,j,k) \text{ (i.e., } \overline{y} = 231\text{)}
    \label{eq:2}
  \end{align}
  then $w \geq z = u\circ (i\ k)$ (i.e., $\overline{z} = 321$).
\end{lem}
\begin{proof}
  Notice that $\dzw(p,\cdot) = \dxw(p,\cdot)$ for $p < j$ and
  $\dzw(p,\cdot) = \dyw(p,\cdot)$ for $p \geq j$.  By
  \lemref{lem:bdf}, $v \leq w$ if and only if $\dvw \geq 0$.  Since
  $x,y \leq w$, $\dxw,\ \dyw \geq 0$.  Combining this with our first
  observation implies that $z \leq w$.
\end{proof}

Now we are ready to prove \thmref{thm:phiinj}. 

\begin{proof}
  First we show that $\pht(\ryw) \subset \rp{yt}{w}$ --- i.e., $yt <
  yt\cdot\pht(t') \leq w$ for all $t'\in\ryw$.  It is clear from the
  definition that $yt \leq yt\cdot\pht(t')$.  So this amounts to showing
  that one of the two dotted arrows in Figure 10 (corresponding to
  relation under the Bruhat order) exists.
  \begin{figure}[htbp]
   \begin{center}
     \xymatrix@dr@C=1pc@W=3.5in{w\ar[d]\ar@{.>}[dr]^? \\
       {yt'}\ar[d]\ar@{.>}[r]^? &
       {yt\pht(t')}\ar[d]^{\pht(t')}\\
       {y}\ar[r]_t\ar[u]^{t'} & {yt}}
     \caption{}
   \end{center}  
  \leavevmode
  \end{figure}\label{fig:genphi}
  
  First suppose $t't = tt'$ --- hence $\pht(t') = t'$.  We wish to
  show that $yt' > yt\pht(t') = ytt' = yt't$.  Letting $t = \tab$
  and $t' = \tcd$, this reduces to showing that $yt'(a) > yt'(b)$.
  Now $t't = tt'$ implies that $\{a,b\}\cap \{c,d\} = \emptyset$.  So
  $yt'(a) > yt'(b)$ if $y(a) > y(b)$.  But this last inequality holds
  by choice of $t$.
  
  Now suppose that $t't \neq tt'$.  In all of the cases in Table
  \ref{tab:dfn} except A.iii and B.iii, $\overline{yt\pht(t')} \leq
  \overline{yt'}$.  Since $yt\pht(t')$ and $yt'$ agree at all indices
  except $a,b,c$, we can use \lemref{lem:comb} to conclude that
  $yt\pht(t') \leq yt' \leq w$.  So $\pht(t')\in \rytw$ for all cases
  except possibly A.iii and B.iii.
  
  In Case A.iii, we know that both $t_{a,c}$, $t_{b,c}\in \ryw$ ---
  i.e., $y < yt_{a,c}, yt_{b,c}\leq w$.  Since $\overline{yt_{a,c}} = 312$
  and $\overline{yt_{b,c}} = 231$, we can therefore invoke
  \lemref{lem:z} to conclude that $w \geq \widehat{321} =
  yt\pht(t')$.  Case B.iii is similar.  
  
  It is clear from Table \ref{tab:dfn} that $\pht(t')$ can share
  at most one index with $t$.  As we have already shown the inclusion
  $\pht(\ryw) \subset \rytw$, we conclude that $\pht(\ryw) \subset
  \rytw \setminus \{t\}$.
  
  Now we show that $\pht$ is an injection.  Suppose $t = t_{\{i,j\}}$,
  $t' = t_{\{j,k\}} \in \ryw$.  No matter which case of Table
  \ref{tab:dfn} we are in, we see that $\pht(t') = t_{\{j,k\}}$ or
  $t_{\{i,k\}}$.  In other words, the index $t'$ doesn't share with
  $t$ must be an index of $\pht(t')$.  In particular, if
  $t',t''\in\ryw$ such that $\pht(t') = \pht(t'')$, the index they
  don't share with $t$ must be the same.  It is then easy to check by
  inspection of Table \ref{tab:dfn} that $\pht(t') = \pht(t'')$
  implies that $t' = t''$.  Since $\pht(t') = t'$ whenever $t$ and
  $t'$ don't share any indices, we conclude that $\pht$ is, in fact,
  an injection.
\end{proof}

The above theorem has as a simple corollary a special case of
Deodhar's conjecture \cite{Deodhar85}.  
\begin{cor}
  $\#\ryw \geq l(w) - l(y)$.
\end{cor}
Various forms and generalizations of the preceding corollary have been
proven by Dyer \cite{Dyer} (arbitrary Coxeter systems), Deodhar
\cite{Deodhar85} (type A), Polo \cite{polo94} (finite Weyl groups) and
Carrell-Peterson \cite{carr94} (crystallographic groups).

Recall from \defref{dfn:tilde} that 
\begin{equation}
    \delxw = \{i, 1 \leq i \leq n: \exists\ j, 1\leq j \leq n,
    \text{ with } t_{\{i,j\}}\in\rxw\}.
  \end{equation}
From the proof of \thmref{thm:phiinj}, we obtain the following:
\begin{cor}\label{cor:delxy}
  For $t\in\trans$, $yt < y \leq w$ implies that $\delyw \subseteq
  \delytw$.
\end{cor}

Next, we see that when $t$ is an adjacent transposition with $wt < w$,
$\pht^{y,w}$ surjects onto $\rytw \setminus \{t\}$.  It would be
interesting to classify all such $t$ for which this happens.

The following fact about the Bruhat order will in useful in the proof
of the next proposition and later in the paper.  An analogous
left-handed version exists.
\begin{lem}[\cite{Hum},Lemma 7.4]\label{lem:hum}
  If $s\in\simref$ and $ws < w$, then $xs \leq w \Leftrightarrow x
  \leq w$.
\end{lem}

\begin{prop}\label{prop:phimax} 
  Let $s_i\in \simref$ (i.e., $s_i$ is an adjacent transposition) and suppose $y,w$
  satisfy the relations $ys_i < y \leq w$ with $ws_i < w$.  Then $\rp{ys_i}{w}
  = \ph_{s_i}(\ryw) \cup \{s_i\}$ and $e_y$ is smooth in $X_w$ if and
  only if $e_{ys_i}$ is smooth in $X_w$.
\end{prop}

\begin{proof}
  The last statement follows immediately from the decomposition of
  $\rp{ys_i}{w}$.  So, taking advantage of \thmref{thm:phiinj}, we
  need only show that if $\tjk\in\rysiw \setminus \{s_i\}$, then
  $\ph_{s_i}^{-1}(t') \neq \emptyset$.
  
  First consider the case where $\{i,i+1\} \cap \{j,k\} = \emptyset$.
  Then $y\tjk = ys_is_i\tjk = ys_i\tjk s_i$.  As $\tjk\in\rysiw$, we can apply
  \lemref{lem:hum} to conclude that $\tjk\in\ryw$.  Finally, since
  $\{i,i+1\} \cap \{j,k\} = \emptyset$, $\ph_{s_i}(\tjk) = \tjk$.  So
  $\ph_{s_i}^{-1}(\tjk) \neq \emptyset$ in this case.
  
  Now, as in the proof of injectivity, we can restrict to the cases
  where $\{i,i+1\} \cap \{j,k\} \neq \emptyset$.  We group into cases
  according to $\overline{y}$, $\overline{ys_i}$ and $t'$.  Cases
  $I,III,V$ have $i \in \{j,k\}$; Cases $II,IV,VI$ have $i+1 \in
  \{j,k\}$.  In each case we show that $t'\in\im\ph_{s_i}$.

  \begin{center}
    \setlength{\extrarowheight}{4pt}
    \begin{tabularx}{5in}{|l|>{$}c<{$}>{$}c<{$}>{$}c<{$}>{$}c<{$}|X|}\hline
      Case & \overline{ys} & \overline{y} & t' & \overline{ys_it'} & $\ph_{s_i}^{-1}(t')$\\\hline
      $I$ & 123 & 213 & t_{a,c} & 321 & Since $w \geq \widehat{321}$,
      $w \geq \widehat{312} = yt_{a,c}$ and $w \geq \widehat{231} =
      yt_{b,c}$.  So we are in Case A.iii and we see that $t_{a,c}
      = \ph_{s_i}(t_{b,c})$.\\\hline
      
      $II$ & 123 & 132 & t_{a,c} & 321 & Case B.iii --- analogous to I.\\\hline

      $III$ & 123 & 213 & t_{b,c} & 132 & By \lemref{lem:hum},
      $\widehat{312} \leq w$ if and only if $\widehat{132} \leq w$.  The
      latter inequality is true since $\overline{ys_it'} = 132$.  So
      $t_{a,c}\in \ryw$.  We're in Case A.i and $t_{b,c} =
      \ph_{s_i}(\tac)$.\\\hline  
      
      $IV$ & 123 & 132 & t_{a,b} & 213 & Case B.i --- analogous to III.\\\hline
  
      $V$& 132 & 312 & t_{a,c} & 231 & $w \geq y = \widehat{312}, w
      \geq ys_it' =
      \widehat{231}$.  Hence, by \lemref{lem:z}, $w \geq \widehat{321}$.
      So $t_{b,c}\in \ryw$.  We're in Case C.i and $t_{a,c} =
      \ph_{s_i}(t_{b,c})$.\\\hline
  
      $VI$ & 213 & 231 & t_{a,c} & 312 & Case D.i --- analogous to V.\\\hline
    \end{tabularx}
  \end{center}
 
  Cases I-VI are the only ones possible since $\overline{ys_i} =
  231,312,321$ cannot have $s,t'$ as hypothesized.
\end{proof}

\begin{prop}\label{prop:phimaxl} 
  Let $s_i\in \simref$ (i.e., $s_i$ is an adjacent transposition) and suppose $y,w$
  satisfy the relations $s_iy < y \leq w$ with $s_iw < w$.  Then $\rp{s_iy}{w}
  = y^{-1}s_i \ph_{s_i}(\ryiwi) s_iy \cup \{y^{-1}s_iy\}$.
\end{prop}
  
\begin{proof}
  By \propref{prop:phimax}, $\rp{y^{-1}s_i}{w^{-1}} =
  \ph_{s_i}(\ryiwi) \cup \{s_i\}$.
  The result follows from the identity $\rxw = x^{-1}\rxiwi x$.
\end{proof}

\begin{cor}\label{cor:singext}  
  If $x \leq w$ is an \msp, $s,s'$ any simple reflections, then $sw <
  w$ or $ws' < w$ imply, respectively, that $sx < x$ or $xs' < x$.
\end{cor}

The preceding corollary is well-known \cite[8.2.10]{BLak}.  However,
\propref{prop:phimax} gives a different proof.

There is one more fundamental property of $\pht$ that we will need to
know for the rest of the paper.  Namely, for a pair of reflections
$t,t'$ where $t\in\im\phtp$, it will be useful to know what we can
say about the membership of $t'$ in $\im\pht$.


\begin{prop}[Reciprocity]\label{prop:reci}
  If $t,t'\in\rxw$, $t \neq t'$, with $l(xt) = l(xt') = l(x) + 1$, then
  $t'\in\im \pht^{xt,w} \Leftrightarrow t\in\im \phtp^{xt',w}$.
\end{prop}
\begin{remark}
  Reciprocity does not necessarily hold if $l(xt') > l(x)+1$.  For
  example, take $x = [1,2,3]$, $w = [3,2,1]$, $t = t_{1,2}$ and $t' =
  t_{1,3}$.  Then $t'\in\im\pht^{xt,w}$ but $t\not\in\im\phtp^{xt',w}$.
\end{remark}
\begin{proof}
  Suppose $t\in\im\phtp$.  We will show that $t'\in\im\pht$.  
  
  First, consider the case where $tt' = t't$.  From the definition of
  $\ph$, we see that $\phtp^{-1}(t) = t$.  So $w \geq xt't = xtt'$.
  This implies that $t' \in \rp{xt}{w}$ and therefore $\pht(t') = t'$.
  
  Now we suppose $tt' \neq t't$.  So $a < b < c$ are determined
  such that $d\not\in \{a,b,c\}$ implies $x(d) = xt(d) = xt'(d)$.
  Let $\overline{x} = \fl_{abc}(x)$.  Note that:
  \begin{enumerate}
  \item By hypothesis, $l(xt) = l(xt') = l(x) + 1$.
  \item If $\overline{x} \in \{231,312,321\}$ then at most one of
    $\tab,\tac,\tbc \in \rxw$; not two.
  \end{enumerate}
  Hence, the cases below are the only ones we need consider.

  \begin{enumerate}
  \item $\overline{x} = 132$.\\
    Then $\{t,t'\} = \{\tab,\tac\}$ and $w \geq
    \widehat{312},\widehat{231}$.  By \lemref{lem:z}, $w \geq
    \widehat{321}$.  So $\tab\in\rp{x\tac}{w}$ and, as we are in Case
    D.ii of \defref{def:phi}, $\tab = \ph_{\tac}(\tab)$.  Similarly,
    $\tbc\in\rp{x\tab}{w}$ and, as we are in Case C.i of
    \defref{def:phi}, $\tac = \ph_{\tab}(\tbc)$.

  \item $\overline{x} = 213$.\\
    The argument is parallel to that in the previous case.

  \item $\overline{x} = 123$.\\
    Here $\{t,t'\} = \{\tab,\tbc\}$.  (Note that $\tac$ is not
    considered since $l(x\tac) > l(x) + 1$.)  Suppose
    $\tab\in\im\ph_{\tbc}$.  From Case B of \defref{def:phi}, this implies
    that
    \begin{equation*}
      \left\{
      \begin{matrix}
        w \geq \widehat{312}\\
        \text{ or }\\
        w \geq \widehat{231}
      \end{matrix}
      \right\}
      \Leftrightarrow
      \left\{
      \begin{matrix}
        \tac\in\rp{x\tab}{w}\\
        \text{ or }\\
        \tbc\in\rp{x\tab}{w}
      \end{matrix}
      \right\}.
    \end{equation*}
    Then, from Case A of \defref{def:phi}, we see that
    $\tbc\in\im\ph_{\tab}$.  The argument is analogous if we instead
    assume $\tbc\in\im\ph_{\tab}$.
  \end{enumerate}
\end{proof}

\section{Preparatory lemmas}
\label{sec:prep}

Let $xt < x \leq w$.  We make the following observation (see
\thmref{fact:basic}):
If 
\begin{equation}
  \#\pht(\rxw) < \#\rxtw - l(xt) + l(x),    
\end{equation}
then $e_{xt}$ is a singular point of $X_w$.
  
The above fact is most conveniently expressed in terms of the
following notation:
\begin{dfn}
  For $x < w$ and $t\in\rxw$, let 
  \begin{equation}
    \qptxw = \rxw \setminus \left(\{t\} \cup \pht(\rp{xt}{w})\right)
  \end{equation}
  denote the set of ``extra'' reflections corresponding to $x$ and $t$.
  We often write $\qabxw$ for $\qptabxw$.
\end{dfn}
If $t'\in\qptxw$, then we say that $t$ and $t'$ are
\textit{incompatible} edges (in the Bruhat graph).  The elements of
$\qptxw$ are ``extra'' edges in the sense that they correspond to an
increase in the dimension of the Zariski tangent space.

The utility of $\qptxw$ is embodied in the following two facts.
\begin{fact}\label{fact:sing}
  If $t,t'\in\rxw$ with $t'\in\qptxw$ and $l(xt) = l(x) + 1$, then
  $x < w$ is singular.
\end{fact}

\begin{fact}\label{fact:enone}
  $x$ is an \msp\ for $w$ if and only if, for every $t\in\rxw$ with
  $l(xt) = l(x) + 1$, $\qptxw \neq \emptyset$.
\end{fact}
Note, however, that if $x$ is a singular point, but not an \msp, then it is
possible that $\qptxw = \emptyset$.  An example is afforded by $w =
[4,2,3,1]$, $x = [1,2,3,4]$ and $t = t_{1,2}$.  Conversely, if $l(xt)
> l(x) + 1$, then we may have $\qptxw \neq \emptyset$ even if $X_w$ is
entirely smooth.  Take, for example, $w = [3,2,1]$, $x = [1,2,3]$, and
$t = t_{1,3}$.

There will be numerous instances in the remainder of the paper where
we do the following:
\begin{enumerate}\label{enum:tech}
\item Assume we have an \msp\ $x$ for $X_w$.
\item Construct some $y = xt'' > x$.
\item Conclude that $y < w$ from the fact that $\area_{t''}(x)$ is shaded.
\item Find incompatible edges $t,t'$ as in \factref{fact:sing} to conclude that $y$ is
  also a singular point of $X_w$.
\item Obtain a contradiction with our first assumption.
\end{enumerate}

The previous technique will allow us to significantly pare down the
possibilities for what $\xti$ looks like for $x$ an \msp.  The
following lemma is the first example of this strategy.
\begin{lem}[Ell Lemma]\label{lem:ell} 
  Let $x \leq w$ and $1 \leq i < j < k \leq n$.
  \begin{enumerate}
   \item If $\fl_{ijk}(x) = 213$ and $\tik,\tjk\in\rxw$, then
     $\tik\in\im \ph_{\tjk}^{x\tjk,w}$ and $\tjk\in\im \ph_{\tik}^{x\tik,w}$.
     
   \item If $\fl_{ijk}(x) = 132$ and $\tij,\tik\in\rxw$, then
     $\tij\in\im \ph_{\tik}^{x\tik,w}$ and $\tik\in\im \ph_{\tij}^{x\tij,w}$.
  \end{enumerate}
\end{lem}

\begin{proof}
  We only prove 1 as the proof for 2 is entirely analogous.\\
  Diagrams for $x$, $x\tik$ and $x\tjk$ are given in
  \figref{fig:ell}.  
  \myfig{ell}{}
  
  We see that $\tik\in\rxw$ implies $w \geq \widehat{312}$ and
  $\tjk\in\rxw$ implies $w\geq \widehat{231}$.  So, by \lemref{lem:z},
  $x\tik\tjk = x\tjk\tij = \widehat{321} \leq w$.  Equivalently,
  $\tjk\in\rp{x\tik}{w}$ and $\tij\in\rp{x\tjk}{w}$.  So, (Case C.ii
  of \defref{def:phi}) $\ph_{\tik}^{x\tik,w}(\tjk) = \tjk$ and (Case D.i of
  \defref{def:phi}) $\ph_{\tjk}^{x\tjk,w}(\tij) = \tik$.
\end{proof}

The next lemma is used frequently.  It gives us criteria for
determining when two reflections are, in fact, incompatible.
\begin{lem}
  Let $\tab \in \rxw$.
  \begin{enumerate}
  \item Patch Incompatibility.\\
    If $\tcd\in\rxw$ with $\{a,b\} \cap\{c,d\} = \emptyset$, then
    $\tab\in\qp{c}{d}$ if and only if $\min(\dxw|_A) = 1$ (with region
    A as in \figref{fig:incom}.1).
  \item Link Incompatibility.\\
    If $\tbc\in\rxw$, then 
    \begin{equation}
      \begin{split}
        \tab\in\qbcxw &\Longleftrightarrow \tbc\in\qabxw\\ 
        &\Longleftrightarrow \min(\dxw|_B) = \min(\dxw|_C) = 0,      
      \end{split}
    \end{equation}
    (where regions B and C are as in \figref{fig:incom}.2).
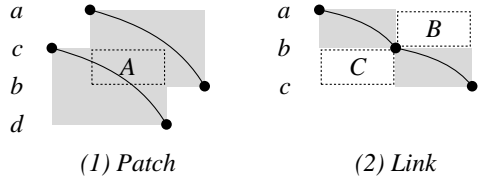
    \myfig{incom}{In 1, we display only one possible configuration where
    $\tab$ and $\tcd$ are patch incompatible.}
  \end{enumerate}
\end{lem}

\begin{proof}
  The proof of Patch Incompatibility is clear.
  To prove Link Incompatibility, it suffices to consider Cases A and B
  of \defref{def:phi}.
\end{proof}

In \factref{fact:sing} we give a sufficient condition for $x$ to be a
singular point of $X_w$ that is expressed in terms of the map $\ph_t$.  Namely,
$x$ is singular if $\qptxw$ is non-empty for some $t\in\rxw$ with
$l(xt) = l(x) + 1$.  As the lemma below shows, many elements of $\rxw$
aren't even candidates to be elements of $\qptxw$.

\begin{lem}\label{lem:disj}
  Let $\tab,\tcd\in\rxw$.  If $\overline{\area_{a,b}} \cap
  \overline{\area_{c,d}} = \emptyset$, then $\tab\not\in\qcdxw$ (i.e.
  --- $\tab\in\im\ph_{\tcd}$).
\end{lem}
\begin{proof}
  For any point $\diamond\in \area_{a,b}(x)$, $\dpw{x\tab}(\diamond) =
  \dxw(\diamond) - 1$.  Similarly for the pair $\tcd$ and $\area_{c,d}(x)$.
  Now, $\overline{\area_{a,b}} \cap \overline{\area_{c,d}} =
  \emptyset$ implies that $\area_{a,b}(x) = \area_{a,b}(x\tcd)$ and
  $\dpw{x\tcd} \geq 1$ on $\area_{a,b}(x\tcd)$.  This implies
  that $\tab\in\rp{x\tcd}{w}$ and $\tab\in\im\ph_{\tcd}$.
\end{proof}

The following lemma is technically useful for \propref{prop:msplo}.

\begin{lem}\label{lem:contain}
  Let $x < w$, $\tab\in\rxw$, $l(x\tab) = l(x) + 1$, and
  $\tcd\in\qabxw$.  If $\pt_x(a)\in\area_{c,d}(x)$ or
  $\pt_x(b)\in\area_{c,d}(x)$, then $x$ is not an \msp\ for $w$.
\end{lem}
\begin{proof}
  First consider the case where both $\pt_x(a),\pt_x(b)\in\area_{c,d}(x)$.
  \myfig{con1}{}
  Suppose there is a point $\pt_x(f)$ in region A of
  \figref{fig:con1}.1.  Choose such an $f$ as small as possible.  Then
  we see that $t_{c,b}$ and $t_{a,f}$ are patch incompatible
  reflections for $x' = xt_{b,d}t_{f,d} \leq w$ and $l(x't_{a,f}) =
  l(x') + 1$.  By \factref{fact:sing}, $x'$ is then singular.  This
  contradicts the fact that $x$ is an \msp\ for $w$.
  
  Now suppose region A of \figref{fig:con1}.1 is empty --- this is
  shown in \figref{fig:con1}.2.  Then $t_{c,b}$ and $t_{a,d}$ are
  incompatible reflections for $xt_{b,d} \leq w$ and
  $l(xt_{b,d}t_{a,d}) = l(xt_{b,d}) + 1$.  Since $x < xt_{b,d}$, this
  contradicts the fact that $x$ is an \msp\ for $w$.

  We now argue the case of $\pt_x(b)\in\area_{c,d}(x)\not\ni\pt_x(a)$.
  (The arguments for $\pt_x(b)\not\in\area_{c,d}(x)\ni\pt_x(a)$ are
  parallel.)

  Clearly $d > b$ and $x(d) > x(b)$.  There are four possibilities
  with regard to the position of $\pt_x(c)$.
  \begin{enumerate}
  \item $c = a$.\\
    We are in Case A.iii of the definition of $\ph$.  Hence,
    $\tcd\in\im\ph_{\tab}$, which is a contradiction.  So this case
    cannot occur.
    
   \item $c > a$, $x(c) > x(a)$.\\
     This case cannot occur as it violates $l(x\tab) = l(x) + 1$.

   \item $c > a$, $x(c) < x(a)$.\\
     This case is depicted in \figref{fig:con2}.1.  Suppose
     $l(xt_{b,d}t_{a,d}) = l(xt_{b,d}) + 1$.  Then $t_{a,d}$ and
     $t_{c,b}$ are patch incompatible for $xt_{b,d} \leq w$.  This
     contradicts the fact that $x$ is an \msp\ for $w$.  If
     $l(xt_{b,d}t_{a,d}) > l(xt_{b,d}) + 1$, then we can argue as in
     \figref{fig:con1}.1 to obtain our contradiction.

   \item $c < a$, $x(c) > x(a)$.\\
     See \figref{fig:con2}.2.  This is analogous to the previous case.
     \myfig{con2}{}

  \end{enumerate}
\end{proof}

\propref{prop:msplo} below gives us our first non-trivial restriction
regarding the composition of $\rxw$.  This proposition will greatly
reduce the amount of work we need to do later on to determine
possibilities for $\xti$.

\begin{prop}\label{prop:msplo}
  Let $x < w$ be an \msp.  If $t\in\rxw$ then $l(xt) = l(x) + 1$.
\end{prop}
\begin{proof}
  Suppose that $t\in\rxw$ and $l(xt) > l(x) + 1$.  We will obtain a
  contradiction.  
  
  Let $t = t_{a,c}$.  Choose $b$ as large as possible such that
  $\pt_x(b)\in\area_{a,c}(x)$.  Note that $\tab,\tbc\in\rxw$ and
  $l(x\tbc) = l(x) + 1$.  Since $x$ is an \msp, we can invoke
  \factref{fact:enone} to find a $t_{e,f}\in\qp{b}{c}$.
  
  Suppose $\overline{\area_{a,b}(x)} \cap \overline{\area_{e,f}(x)} =
  \emptyset$.  Since $t_{a,c}\in\rxw$, $\area_{a,c}$ is shaded so
  $\tbc\in \rp{x\tab}{w}$.  Hence $t_{b,c}$ and $t_{e,f}$ are
  incompatible for $x\tab \leq w$ and $l(x\tab\tbc) = l(x\tab) + 1$.
  This contradicts $x\in\maxsing$.
  
  Otherwise, $\overline{\area_{e,f}}$ overlaps both
  $\overline{\area_{a,b}}$ and $\overline{\area_{b,c}}$, so, by
  \lemref{lem:contain}, we are in one of the following two scenarios.
  \begin{enumerate}
  \item $e = b$.\\
    By choice of $b$, $f\not\in\area_{b,c}(x)$.  (Note that $f \neq
    c$.)  So either $f > c$, $x(f) < x(c)$ or $f < c$, $x(f) > x(c)$
    (the latter case is shown in \figref{fig:lone}.1.).  In either
    case, we can apply the Ell Lemma~\ref{lem:ell} to conclude that
    $t_{e,f}\in\im\ph_{\tbc}$.  This contradicts the choice of
    $t_{e,f}$.  \myfig{lone}{}
  \item $f = b$.\\
    Since $\tac\in\rxw$, for $t_{e,f}$ to be an element of
    $\qp{b}{c}$, we need $e < a$, $x(e) < x(a)$ and $\dxw = 0$ for
    some point in each of regions A and B in \figref{fig:lone}.2.  But
    then $t_{e,b}$ and $t_{b,c}$ are link incompatible for $x\tab
    \leq w$.  Furthermore, by having chosen $b$ as large as possible, we
    ensure that $l(x\tab\tbc) = l(x\tab) + 1$.  This contradicts
    $x\in\maxsing$.
  \end{enumerate}
\end{proof}

We are now able to give a graphical description of all possible pairs
$\tab,\tcd\in\rxw$ such that $\tcd\in\qabxw$ when $x\in\maxsing$.  The
result follows immediately from \lemref{lem:disj} and
\propref{prop:msplo} along with the Ell Lemma~\ref{lem:ell}.  Note that by
\propref{prop:reci}, we can assume, without loss of generality, that
$a < c$.
\begin{cor}\label{cor:grap}
  If $\tab\in\rxw$ and $\tcd\in\qabxw$, then the relative positions of
  $\area_{a,b}$ and $\area_{c,d}$ are one of the ones shown in
  \figref{fig:relpos}.
  \myfig{relpos}{}
\end{cor}

This greatly simplifies our future investigations.  We now use
\corref{cor:grap} and \lemref{lem:rhombus} to prove the following
crucial lemma.
\begin{lem}[Cross Lemma]\label{lem:cross} 
  Let $x < w$ be an \msp\ and suppose $1 \leq i < j < k < l \leq
  n$ such that $\fl_{ijkl}(x) = 2143$.  If $t_{j,k}\in\rxw$
  and $t_{i,l}\in\qp{j}{k}$, then $t_{i,k},t_{j,l}\in\rxw$.
\end{lem}
\begin{proof}
  We can visualize the situation as in \figref{fig:cross}.1.
  \myfig{cross}{} 
  Since $t_{i,l}\in\qp{j}{k}$, there is necessarily a point
  $\triangle$ in region A for which $\dxw(\triangle) = 1$.  Suppose
  $t_{i,k}\not\in\rxw$.  Then there is a point $\square$ in region B
  such that $\dxw(\square) = 0$.  Then we can apply
  \lemref{lem:rhombus} (with $\alpha,\beta \geq 1$, $\gamma = 1$) and
  \propref{prop:msplo} to conclude that there is a point $\pt_x(p)$ of
  $\digx$ in region $C$ (see \figref{fig:cross}.2).  If we choose
  $\square$ to be as low as possible in our diagram, then $\dxw|_D
  \geq 1$ (see \figref{fig:cross}.3).  But then $t_{i,l}$ and
  $t_{j,k}$ are patch incompatible for $xt_{p,k} \leq w$.  This
  contradicts $x\in\maxsing$.  Therefore $\tik\in\rxw$ and we can
  shade the entire region B.

  To shade the lower left corner, apply the preceding argument to
  $x^{-1}$ and $w^{-1}$.
\end{proof}

\section{Restrictions on $\xti$}
\label{sec:rest}

Recall that $\xti$ and $\wti$ are the restrictions of $x$ and $w$ to
those positions in $\delxw$ (see \defref{dfn:tilde}).  In order to
determine the structure of $\maxsing$, we first prove the following
necessary conditions on $\xti$ for any \msp\ $x$ for $w$.

\begin{thm}\label{thm:l7}
  If $x < w$ is an \msp, then $\xti$ is 231- and 312-avoiding.
\end{thm}
\begin{thm}\label{thm:stair}
  If $x < w$ be an \msp, then $\xti$ is 1234-avoiding.
\end{thm}

These two theorems are, in fact, almost
enough to describe $\xti$ for any \msp\ $x$.


\subsection{Technical lemmas regarding $\xti$}
\label{subsec:techl7}

In order to streamline the proof of \thmref{thm:l7}, we first present
two technical lemmas.  These lemmas simply show that $\digxt$, for $x$
an \msp, must avoid certains patterns of points and shading.

\begin{lem}\label{lem:dagger}
  If $x < w$ is an \msp, then $\digxt$ does not contain either of the
  configurations in \figref{fig:dagger0} (regardless of whether or not
  these reflections are incompatible).
  \myfig{dagger0}{}
\end{lem}
\begin{proof}
  We only prove that $\digxt$ must avoid the configuration in
  \figref{fig:dagger0}.1.  The proof of the other case is parallel.
  
  First suppose that $\min(\dxw|_A) = 1$ (see
  \figref{fig:dagger0}.1).  Then by the Cross Lemma~\ref{lem:cross}, $t_{a,\gamma} \in
  \rxw$.  But this contradicts \propref{prop:msplo} as
  $l(xt_{a,\gamma}) > l(x) + 1$.  We get a similar contradiction if
  $\min(\dxw|_B) = 1$.  So we can henceforce assume that our
  configuration is actually as in \figref{fig:dagger1} --- i.e.,
  $\min(\dxw|_A)$, $\min(\dxw|_B) \geq 2$.
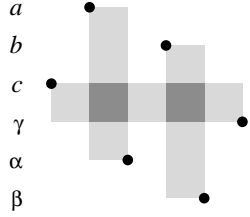
  \myfig{dagger1}{Recall that dark shading denotes regions on which
  $\dxw \geq 2$.}
  
  Since $x$ is an \msp, we must have some reflection
  $t_{d,\delta}\in\qp{c}{\gamma}$.  Suppose
  \begin{align}\label{eq:3}
    \overline{\area_{d,\delta}} \cap
    \overline{\area_{a,\alpha}} &= \emptyset\text{ or }\\
    \overline{\area_{d,\delta}} \cap \overline{\area_{b,\beta}} &=
    \emptyset.\label{eq:6}
  \end{align}
  In the case of (\ref{eq:3}), $t_{c,\gamma}$ and $t_{d,\delta}$ are
  incompatible for $xt_{a,\alpha} \leq w$.  This contradicts
  $x\in\maxsing$.  (Similarly for (\ref{eq:6}).)  So these
  intersections must be non-empty.  It is clear from
  \figref{fig:dagger1} that for these intersections to be non-empty,
  we need $t_{d,\delta}$ to be patch (rather than link) incompatible
  with each of $t_{c,\gamma}$, $t_{a,\alpha}$ and $t_{b,\beta}$.  By
  \propref{prop:msplo} and because the intersections in (\ref{eq:3})
  and (\ref{eq:6}) must be non-empty, it is readily seen that we
  require $\delta < \alpha$ and $d > b$.  Hence, there
  are only four possible ways in which $\area_{d,\delta}$ may overlap
  $\area_{a,\alpha}$, $\area_{b,\beta}$ and $\area_{c,\gamma}$.  These
  are shown in \figref{fig:dagger2}.  \myfig{dagger2}{}
  \begin{enumerate}
  \item $d > c$, $\delta > \gamma$.\\
    If $\min(\dxw|_C) \geq 2$ then $t_{c,\gamma}$ and $t_{d,\delta}$
    are patch incompatible for $xt_{a,\alpha} \leq w$.  This
    contradicts $x\in\maxsing$.  Otherwise, we can apply the Cross Lemma~\ref{lem:cross}
    to $t_{d,\delta}$ and $t_{a,\alpha}$ to conclude that
    $t_{a,\delta}\in\rxw$.  This contradicts \propref{prop:msplo}.

  \item $d < c$, $\delta < \gamma$.\\
    The argument is parallel to the previous case.
    
  \item $d < c$, $\delta > \gamma$.\\
    If $\min(\dxw|_C) = 1$ or $\min(\dxw|_D) = 1$ then we can apply
    the Cross Lemma~\ref{lem:cross} to $t_{d,\delta}$ and $t_{a,\alpha}$ to conclude
    that $t_{a,\delta}\in\rxw$.  This contradicts
    \propref{prop:msplo}.  The only alternative is that $\dxw \geq 2$
    on regions C and D.  But then $t_{d,\delta}$ and $t_{c,\gamma}$
    are patch incompatible for $xt_{a,\alpha} \leq w$.  This
    contradicts $x\in\maxsing$.

  \item $d > c$, $\delta < \gamma$.\\
    Here $t_{c,\gamma}$ and $t_{d,\delta}$ are patch incompatible
    for $xt_{a,\alpha} \leq w$.  This contradicts $x\in\maxsing$.
  \end{enumerate}
\end{proof}

\begin{lem}\label{lem:xi}
  If $x < w$ is an \msp, then $\digxt$ does not contain the
  configuration in \figref{fig:xi1}.
  \myfig{xi1}{}
\end{lem}

\begin{proof}
  Since $x$ is an \msp\ for $w$, there exists some
  $t_{d,\delta}\in\qp{c}{\gamma}$.  Clearly, if
  \begin{equation}\label{eq:10}
      \overline{\area_{d,\delta}} \cap \overline{\area_{\alpha,a}}
      = \emptyset,
  \end{equation}
  then $t_{c,\gamma}$ and $t_{d,\delta}$ are (patch or link) incompatible
  reflections for $xt_{\alpha,a} \leq w$.  This would contradict
  $x\in\maxsing$.
  
  So, to ensure that (\ref{eq:10}) does not hold, we need $\pt_x(d)$
  in region A of \figref{fig:xi2}.1 and $\pt_x(\delta)$ in region B.
  Here we are including the possibilities that $d = a$ or $\delta =
  c$.  Note that (as is shown in \figref{fig:xi2}.1)
  \propref{prop:msplo} requires that $x(\delta) < x(b)$ and $x(d) >
  x(\beta)$.  Clearly if $d = a$ then $\delta \neq c$ and vice versa.
  Hence, by symmetry, we can treat only the cases where $\delta \neq
  c$.  These two cases are illustrated in Figures \ref{fig:xi2}.2 and
  \ref{fig:xi2}.3.
  \myfig{xi2}{} 
  
  For both cases, we can apply the Cross Lemma~\ref{lem:cross} to $t_{d,\delta}$ and
  $t_{c,\gamma}$ to conclude that $t_{d,\gamma}\in\rxw$.  Then
  $l(xt_{d,\delta}) > l(x) + 1$, which contradicts \propref{prop:msplo}.
\end{proof}

\subsection{Proofs of Theorems \ref{thm:l7} and \ref{thm:stair}}
\label{sec:l7stair}

\begin{proof}[\thmref{thm:l7}]
  By passing to inverses, it is enough to prove that $\xti$ is
  231-avoiding.  So choose $a,b,c\in\delxw$ with $1 \leq a < b < c
  \leq n$ such that $\fl_{abc} = 231$.
  \myfig{L7}{$\tab\in\rxw$.}
  \begin{enumerate}
  \item Assume $\tab \in \rxw$.\label{L7CaseI}\\
    We have the situation of \figref{fig:L7}.  
    
    By definition of $\delxw$, there exists a $d\in\delxw$ with
    $\tcd\in\rxw$.  We'll assume that $c < d$ as all cases where $d <
    c$ are analogous to one of the cases we cover by transposing over
    the anitdiagonal.  Clearly $\overline{\area_{a,b}} \cap
    \overline{\area_{c,d}} = \emptyset$.
    
    Since $x$ is an \msp, there exists a $\talbt\in\qabxw$.  We claim that
    \begin{equation}
    \overline{\area_{\alpha,\beta}} \cap
    \overline{\area_{c,d}} \neq \emptyset.
    \end{equation}
    Suppose the intersection is empty.  Then $\tab$ and $\talbt$ are
    (patch or link) incompatible for $x\tcd \leq w$.  This contradicts
    $x\in\maxsing$.  So we may assume that this intersection is
    non-empty.  There are three cases according to whether $\tab$ and
    $\talbt$ are
    \begin{enumerate}
    \item link incompatible
    \item patch incompatible with $\alpha < a$.
    \item patch incompatible with $\alpha > a$.
    \end{enumerate}

    We only describe the arguments explicitly in the case of link
    incompatibility --- the arguments are similar in the latter two
    cases.  We will argue only $b = \alpha$ as the case of $\beta = a$
    is analogous.
    
    By Propositon~\ref{prop:msplo}, there are three possibilities for
    the relative positions of $\pt_x(d)$ and $\pt_x(\beta)$.  They are
    displayed in \figref{fig:joyl}.  \myfig{joyl}{}
    \begin{enumerate}
    \item $d > \beta$ (i.e., $x(d) < x(\beta)$).\\
      We have that $\tab$ and $t_{\alpha,c}$ are link incompatible for
      $xt_{c,\beta}\leq w$.  This contradicts $x\in\maxsing$.
    \item $d = \beta$ (i.e., $x(d) = x(\beta)$).\\
      The argument is the same as in the previous case.
    \item $d < \beta$ (i.e., $x(d) > x(\beta)$).\\
      If $\min(\dxw|_C) \geq 2$, then $\tab$ and $\talbt$ are link
      incompatible for $x\tcd\leq w$.  This contradicts $x\in\maxsing$.
      Otherwise, we can apply the Cross Lemma~\ref{lem:cross} to $\tcd$ and
      $\talbt$ to conclude that $t_{c,\beta}\in\rxw$.  Then, as in
      \figref{fig:joyl}.3, 
      $\tab$ and $t_{b,c}$ are link incompatible for $xt_{c,\beta}
      \leq w$.  This contradicts $x\in\maxsing$.
    \end{enumerate}
    
  \item Assume $\tab \not\in \rxw$.\\
    Since $a,b,c\in\delxw$, we can find $\alpha,\beta,\gamma$ such
    that $t_{\{a,\alpha\}}, t_{\{b,\beta\}}, t_{\{c,\gamma\}}\in\rxw$.
    If \myfig{tto1}{}

    \begin{enumerate}
    \item $\pt_x(\alpha)$ is in region A or A' of \figref{fig:tto1}.1
      or
    \item $\pt_x(\beta)$ is in region C or C' of \figref{fig:tto1}.2
      or
    \item $\pt_x(\gamma)$ is in region E or E' of \figref{fig:tto1}.3,
    \end{enumerate}
    then we can reduce to the previous case (of $\tab\in\rxw$) or we
    violate \propref{prop:msplo}.  So we must have
    \begin{enumerate}
    \item $\pt_x(\alpha)\in B \cup B'$
    \item $\pt_x(\beta)\in D \cup D'$
    \item $\pt_x(\gamma)\in F \cup F'.$
    \end{enumerate}
    As the argument of $\pt_x(\gamma) \in F$ is analogous, we will
    assume $\pt_x(\gamma) \in F'$.
    
    We now argue the four cases we have left according to whether
    $\pt_x(\alpha) \in B$ and whether $\pt_x(\beta) \in D$.
    \myfig{last4}{}
    \begin{enumerate}
    \item $\pt_x(\alpha) \in B$, $\pt_x(\beta) \in D$.\\
      See \figref{fig:last4}.1.  Note that we can reduce
      to Case \ref{L7CaseI} if $x(\alpha) > x(\beta)$.  So we assume
      $x(\alpha) < x(\beta)$.  We have already shown in
      \lemref{lem:xi} that this configuration contradicts $x\in\maxsing$.

    \item $\pt_x(\alpha) \in B$, $\pt_x(\beta) \in D'$.\\
      Since $x$ is an \msp, there exists some
      $t_{d,\delta}\in\qp{b}{\beta}$.  As $\min(\dxw|_A) = 0$, one can
      see in \figref{fig:last4}.2 that $\overline{\area_{d,\delta}}
      \cap \overline{\area_{\alpha,a}} = \emptyset$.  Hence,
      $t_{d,\delta}$ and $t_{b,\beta}$ are (patch or link)
      incompatible for $xt_{\alpha,a} \leq w$.
      
    \item $\pt_x(\alpha) \in B'$, $\pt_x(\beta) \in D$.\\
      See \figref{fig:last4}.3.  We have already shown in
      \lemref{lem:dagger} that this case contradicts $x\in\maxsing$.
      
    \item $\pt_x(\alpha) \in B'$, $\pt_x(\beta) \in D'$.\\
      See \figref{fig:last4}.4.  Note that if $\alpha > \beta$, then
      in this case we can reduce to Case \ref{L7CaseI} (with
      $\fl_{b\beta\alpha}(x) = 231$).  So we assume $\alpha < \beta$.
      We have already shown in \lemref{lem:dagger} that this case
      contradicts $x\in\maxsing$.
    \end{enumerate}
  \end{enumerate}
  This completes the proof that $\xti$ is 231- and 312-avoiding.
\end{proof}

\begin{proof}[\thmref{thm:stair}]
  Suppose we have $a < b < c < d$ with $a,b,c,d\in\delxw$ and
  $\fl_{abcd}(x) = 1234$.  We will obtain a contradiction.

  By \thmref{thm:l7}, no points of $\digxt$ may occur in regions I or
  II of \figref{fig:stair}.1.  
  \myfig{stair}{}
  Since $a,d\in\delxw$, there exist $b',c'$ such that
  $t_{\{a,b'\}},t_{\{c',d\}}\in\rxw$.  As $x$ is an \msp, there also
  exists some $t_{e,f}\in\qp{a}{b'}$.  Using \propref{prop:msplo}, it
  is easy to check that if
  \begin{equation}
    b'\not\in \overline{A \cup B}\text{ or }c'\not\in\overline{C \cup D}\text{ or }e\neq b'\text{ or } f \neq c',
  \end{equation}
  then $t_{e,f}$ and $t_{a,b'}$ are incompatible for $xt_{c',d} \leq
  w$.  This contradicts $x\in\maxsing$.  So, in particular, $t_{e,f} =
  t_{b',c'} \in \rxw$.  A typical (allowable) pair of positions for
  $b'$ and $c'$ is shown in \figref{fig:stair}.2.
  
  As we have now constructed $a < b' < c' < d$ with
  $t_{a,b'},t_{b',c'},t_{c',d}\in\rxw$, we will assume that $b$ and
  $c$ were chosen initially such that $\tab,\ \tbc,\ \tcd\in\rxw$.
  Note that by the above construction, we can assume that
  $\tab\in\qp{b}{c}$.  Suppose $t_{b,c} \not\in\qp{c}{d}$.  Then
  $t_{a,b}$ and ($t_{b,c}$ or $t_{b,d}$) are link
  incompatible for $xt_{c,d} \leq w$.  This contradicts
  $x\in\maxsing$.
  
  So our diagram looks like that pictured in \figref{fig:stair2}.1 and
  we have $t_{a,b} \in \qbcxw, t_{b,c} \in \qcdxw$.  \myfig{stair2}{}
  
  Therefore, we can find a point in each of the regions U \& V such
  that $\dxw = 0$.  Choose the point in region U to be as low as
  possible.  Choose the point in region V to be as far right as
  possible.  Such points are shown in \figref{fig:stair2}.2.  Apply
  \lemref{lem:rhombus} to the rectangle determined by these two points
  with $\alpha,\beta \geq 1$ and $\gamma = 0$.  This, along with
  \propref{prop:msplo}, implies that there is another point $\pt_x(p)$
  in either region P or Q.  Without loss of generality, assume it is
  in region P.  By having chosen the point in region U as low as
  possible, we find that $t_{p,d}\in\rxw$ (see \figref{fig:stair2}.3).
  Hence, $t_{a,b}$ and $t_{b,c}$ are link incompatible for $xt_{p,d}
  \leq w$.  This contradicts $x\in\maxsing$.
\end{proof}

\section{Restrictions on $\wti$}
\label{sec:wtrest}

Combining \thmref{thm:l7} and \thmref{thm:stair},
we see that if $x\in\maxsing$, then
\begin{equation}
  \xti = [k,\ldots,1,k+l,\ldots,k+1,k+l+m,\ldots,k+l+1]
\end{equation}
for some $k,l,m \geq 0$.  If two out of three of $k,l,m$ are $0$, then
$\xti$ is strictly decreasing, so $x\leq w$ implies that $x = w$.  But
then $x$ cannot possibly be an \msp.  So now we determine the
possible values of $k,m$ in \propref{prop:xtab} when $l = 0$ and
the possible values of $k,l,m$ in \propref{prop:xtabc} when $k,l,m >
0$.  In each proposition, we also determine what $\wti$ must be to
allow $\xti$ to be singular.

We know from \propref{prop:rxtwt} that $x\in\maxsing$ if and only if
$\xti\in\maxsingt$.  Hence, for the remainder of this section, we will
only consider the case where $\xti = x$ and $\wti = w$.

\subsection{Two decreasing sequences in $\xti$}
\label{subsec:xtab}

\begin{prop}\label{prop:xtab}
  Let $x\in\maxsing$ with $\xti = x$ and $\wti = w$.  Suppose that $x$
  consists of exactly two decreasing sequences:
  \begin{equation}
    x = [k,\ldots,1,k+m,\ldots,k+1],
  \end{equation}
  for some $k,m\geq 1$.  Then
  \begin{enumerate}
    \item $k,m \geq 2$ \text{ and }
    \item $w = [k+m,k,\ldots,2,k+m-1,\ldots,k+1,1]$ (shown in
      \figref{fig:xtabno2}.3).
  \end{enumerate}
\end{prop}

\begin{proof}
  For brevity in the following, we use the convention that
  $\alpha,a,a'\in[1,\ldots,k]$ and $\beta,b,b'\in[k+1,\ldots,k+m]$.
  
  Condition 1: If $k = 1$ or $m = 1$, then by \lemref{lem:ell},
  $\qpt{1}{k+1} = \emptyset$.  This contradicts \factref{fact:enone}.
  
  Condition 2: We split this proof into proving the following facts:
  \begin{enumerate}
  \item $\tti_{1,k+1},\tti_{k,k+m}\in\rxw$.
  \item $\ttab\in\rxw$ for all $1\leq a\leq k$ and $k+1\leq
    b\leq k+m$.
  \item $\dxw \leq 1$.
  \item $w$ is as in the statement of Condition 2.
  \end{enumerate}

  We now prove these claims.
  \renewcommand{\labelenumi}{Step \theenumi.}
  \begin{enumerate}
  \item $\tti_{1,k+1},\tti_{k,k+m}\in\rxw$.\\
    Assume $\tti_{1,k+1}\not\in\rxw$.  We will obtain a
    contradiction.
    
    By \propref{prop:lblf}, we can find $\alpha,\beta$ such that
    $\tti_{1,\beta},\tti_{\alpha,k+1}\in\rxw$ (see \figref{fig:xtab}).
    Choose $\alpha$ as large, $\beta$ as small as possible subject to
    this restriction.

    \myfig{xtab}{} 
    
    If $\tti_{1,\beta}\in\qpt{\alpha}{k+1}$, then an application of
    the Cross Lemma~\ref{lem:cross} would offer the desired
    contradiction.  So assume that this is not the case (i.e., assume
    $\dxw \geq 2$ on region R of \figref{fig:xtab}).
  
    Since $x$ is an \msp, by \factref{fact:enone}, we can find some
    $\ttab\in\qpt{1}{\beta}$.  Recall that we chose $\alpha$ as large
    as possible such that $\tti_{\alpha,k+1}\in\rxw$.  It follows
    then that $a \leq \alpha$.  Similarly, our choice of $\beta$ as
    small as possible such that $\tti_{1,\beta}\in\rxw$, in
    conjunction with the Cross Lemma~\ref{lem:cross} and Ell
    Lemma~\ref{lem:ell}, implies that $b > \beta$.  Suppose $a =
    \alpha$.  This is depicted in \figref{fig:xtab1}.1.  We see that
    $t_{k+1,\beta}$ and $t_{1,b}$ are patch incompatible for
    $xt_{\alpha,k+1} \leq w$.  This contradicts $x\in\maxsing$.  So we
    may assume $a < \alpha$ as in \figref{fig:xtab1}.2.
    \myfig{xtab1}{}
    
    Suppose that $\dxw \geq 2$ on region A.  Then $\tab$ and
    $t_{1,\beta}$ are patch incompatible for $xt_{\alpha,k+1} \leq w$.
    This contradicts $x\in\maxsing$.  So there is at least one point
    in region A for which $\dxw$ has value 1.  Now we can apply the
    Cross Lemma~\ref{lem:cross} to the patch incompatible pair
    $\tti_{\alpha,k+1}, \ttab$ to conclude that $\dxw \geq 1$ on
    regions B and C.  We display this knowledge in
    \figref{fig:xtab2}.1.
    \myfig{xtab2}{}
  
    Now suppose that there is a point in region D for which $\dxw =
    1$.  Then $t_{k+1,\beta}$ and $t_{1,b}$ are patch incompatible
    reflections for $xt_{\alpha,k+1} \leq w$.  This contradicts
    $x\in\maxsing$.  Since, by construction, $\ttab\in\qpt{1}{\beta}$,
    the only possibility left is that $\min(\dxw|_E) = 1$ (as in
    \figref{fig:xtab2}.2).  We can now apply the Cross Lemma~\ref{lem:cross} to
    $\tti_{1,b}$ and $\tti_{a,k+1}$ to conclude that $\dxw \geq 1$ on region
    F.  Hence $\tti_{1,k+1}\in\rxw$ as claimed.
  
    The proof that $\tti_{k,k+m}\in\rxw$ is entirely analogous when one
    uses $\dxw'$ from (\ref{eq:dualdxw}).

  \item $\tti_{a,b}\in\rxw$ for all $1\leq a\leq k$ and $k+1\leq
    b\leq k+m$.\\
    By the previous step, we know that we can shade regions I and II
    in \figref{fig:xtab2}.3.  For every $a,b$, by the definition of
    $\delxw$, we can shade the corresponding regions U and V,
    respectively.  This completes the claim.

  \item $\dxw \leq 1$.\\
    Suppose, on the contrary, that $\dxw \geq 2$ for some point on
    region A in \figref{fig:xtabno2a}.1.
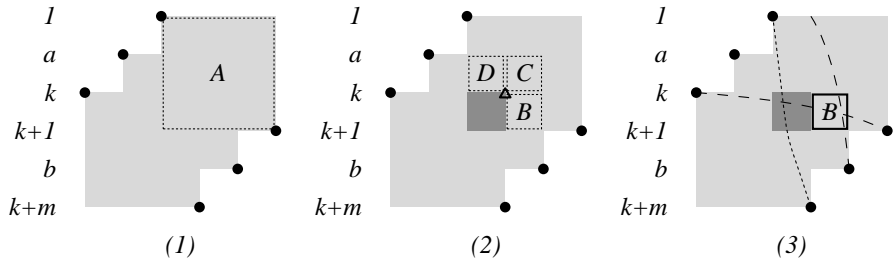
    \myfig{xtabno2a}{We have displayed the case of $a < k$, but the
      argument holds for $a = k$ too.}
    
    Since $\dxw$ is non-decreasing as we move down or left in region
    A, we can assume that $\dxw(\triangle) \geq 2$ for $\triangle =
    (k,k+1)$.  But then there must be some $a,b$
    with $1 < a \leq k$ and $k + 1 < b \leq k+m$ with either $a < k$
    or $b < k+m$ and $\ttab\in\qpt{1}{k+1}$ (see
    \figref{fig:xtabno2a}.2).  Note that $\min(\dxw|_{B\cup C \cup D})
    = 1$ by choice of $\ttab$.  If $\min(\dxw|_B) = 1$ then
    $t_{k,k+1}$ and $t_{1,b}$ are patch incompatible for $xt_{1,k+m}
    \leq w$ (see \figref{fig:xtabno2a}.3).  This contradicts
    $x\in\maxsing$.  So we can assume $\dxw|_B \geq 2$.  Since
    $\min(\dxw|_{B\cup C \cup D}) = 1$, and $\dxw$ is non-decreasing
    in region A as move down or left, we can now assume that
    $\min(\dxw|_C) = 1$.  But then $t_{1,b}$ and $t_{a,k}$ are
    patching incompatible for $xt_{k,k+1}\leq w$ (see
    \figref{fig:xtabno2}.1).  This contradicts $x\in\maxsing$.  So
    $\dxw|_A \leq 1$ as claimed.
    \myfig{xtabno2}{}

  \item By the previous step, there is at most one point of $\digw$
    in region A.  But as $w > x$, $\xti = x$ and $\wti = w$, this
    fixes all the remaining points and we see that
  \begin{equation}
    w = [k+m,k,\ldots,2,k+m-1,\ldots,k+1,1],
  \end{equation}
  as claimed.  This is displayed in \figref{fig:xtabno2}.2.
\end{enumerate}
\end{proof}

\subsection{Three decreasing sequences in $\xti$}
\label{subsec:xtabc}

We repeat the task of the previous section when $\xti$ consists of
three decreasing sequences rather than just two.
\begin{prop}\label{prop:xtabc}
  Let $x\in\maxsing$ with $\xti = x$ and $\wti = w$.  Suppose that $x$
    consists of exactly three decreasing subsequences:
  \begin{equation}
    x = [k,\ldots,1,k+l,\ldots,k+1,k+l+m,\ldots,k+l+1],
  \end{equation}
  for some $k,l,m\geq 1$.  Then
  \begin{enumerate}
  \item $l \geq 2$,
  \item $l = 2$ if $k > 1$ or $m > 1$,\label{item:hmm}
  \item \begin{equation*}
    \begin{split}\label{eq:15}
      w = [k+l,k,\ldots,&2,k+l+m,k+l-1,\ldots,k+2,\\
          &1,k+l+m-1,\ldots,k+l+1,k+1].      
    \end{split}
    \end{equation*}
    (Shown in \figref{fig:wtiabc1}.)
\end{enumerate}
\end{prop}

\begin{proof}
  Again for purposes of brevity, we'll assume throughout this proof
  that $a,a',\alpha\in[1,\ldots,k]$, $b,b',\beta\in[k+1,\ldots,k+l]$ and
  $c,\gamma\in[k+l+1,\ldots,k+l+m]$.  We now prove a series of claims
  elucidating the structure of $\rxw$.

  As the chain of reasoning has several steps, we summarize them here
  before beginning.  
  \begin{enumerate}
  \item There exist $a,b,c$ such that $\ttab,\ttbc\in\rxw$.
  \item $\ttab$ and $\ttbc$ are link incompatible.
  \item Such $a$ and $c$ exist for any $b$.
  \item $\ttab,\ttbc\in\rxw$ for any triple $a,b,c$.
  \item Conclude that we can shade the diagram as in \figref{fig:xtabclen1}.1.
  \item $\dxw \leq 1$ in \figref{fig:xtabclen1}.1.
  \item Condition 1 holds.
  \item $w$ is as claimed in Condition 3.
  \item Condition 2 holds.
  \end{enumerate}
  \renewcommand{\labelenumi}{Step \theenumi.}
  \begin{enumerate}
  \item There exist $a,b,c$ such that $\ttab,\ttbc\in\rxw$.\\
    Suppose there is no such triple of indices.  Then by
    the definition of $\delxw$, for given $a,c$ there exist $b,b'$ ($b\neq
    b'$) such that $\ttab,\tti_{b',c}\in\rxw$.  By the Ell
    Lemma~\ref{lem:ell}, along with the assumptions that $x$ is an
    \msp\ and that such triples do not exist, we can find $\alpha \neq
    a$ and $\beta\neq b'$ such that $\tti_{\alpha,\beta}\in\qabxw$
    (illustrated in \figref{fig:xtabc1}).

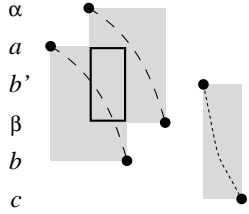
    \myfig{xtabc1}{We have displayed the case of $b' < \beta < b$,
      $\alpha < a$, but there are other possibilities.}
    
    But then $\tab$ and $\talbt$ are patch incompatible for $xt_{b',c}
    \leq w$.  This contradicts $x\in\maxsing$.  So we must be able to
    find a triple as claimed.
    
  \item If $\ttab,\ttbc\in\rxw$ then $\ttab\in\qbcxw$.\label{item:6}\\
    Suppose, on the contrary, that $\ttab\in\im\ph_{\ttbc}$.  
    Then $\dxw \geq 1$ on either all of region A or all of region B in
    \figref{fig:xtabc2a}.1.  
    \myfig{xtabc2a}{} 
    
    Assume that $\dxw \geq 1$ on region A.  Now, since $x$ is an \msp,
    by \factref{fact:enone} there exists $\ttpq\in\qbcxw$.  We now
    consider the two possibilities for the relative positions of
    $\ttpq$ and $\ttbc$.  
    
    Suppose that $\ttpq$ and $\ttbc$ are link incompatible --- i.e., we
    have $q = b$ (\figref{fig:xtabc2a}.2).  For $\tti_{p,q}$ to be link
    incompatible with $\ttbc$, we need $p > a$ (as depicted in
    \figref{fig:xtabc2a}.2) since we are assuming $\min(\dxw|_A) \geq
    1$.  Additionally, as $\ttpq\in\qbcxw$, $\dxw$ must have value 0
    for at least one point on each of regions C and D.  Thus $t_{p,b}$
    and $\tbc$ are link incompatible for $x\tab \leq w$.  This
    contradicts $x\in\maxsing$.
    
    On the other hand, $\ttpq$ and $\ttbc$ may be patch incompatible.
    Then there are four possibilities for the relative positions of
    $\pt_x(p)$, $\pt_x(b)$, $\pt_x(q)$ and $\pt_x(c)$ depending on
    whether $p < b$ and whether $q < c$ (see \figref{fig:xtabc2c}).  In
    each situation, $\tpq$ and $\tbc$ are patch incompatible for
    $x\tab \leq w$.  This contradicts $x\in\maxsing$.  
    \myfig{xtabc2c}{}

    We have obtained a contradiction for every scenario in which $\dxw
    \geq 1$ on region A.  Arguing similarly if $\dxw \geq 1$ on region
    B, we conclude that $\ttab\in\qpt{b}{c}$.

  \item Given $\beta$, there exist $\alpha,\gamma$ such that
    $\tti_{\alpha,\beta},\tti_{\beta,\gamma} \in \rxw$.\label{item:4}\\
    By Step 1, there exist $a,b,c$ such that $\ttab,\ttbc\in\rxw$.
    If $b = \beta$ then we are done --- so assume not.  We can at
    least find a $q$ with $\tti_{\{\beta,q\}}\in\rxw$.  Without loss
    of generality, assume $q = \gamma$ for some $\gamma > \beta$.  We
    split into cases according to whether $\gamma < c$, $\gamma = c$
    or $\gamma > c$.  These are depicted in \figref{fig:xtabc3}.  Note
    that by the previous step, $\ttab\in\qbcxw$, so $\min(\dxw|_B) =
    0$.  In addition, if $\min(\dxw|_A) \geq 1$, then
    $\tti_{a,\beta}\in\rxw$ as desired.  So in the following
    arguments (and \figref{fig:xtabc3}), we assume $\min(\dxw|_A) = 0$
    and derive a contradiction.
    
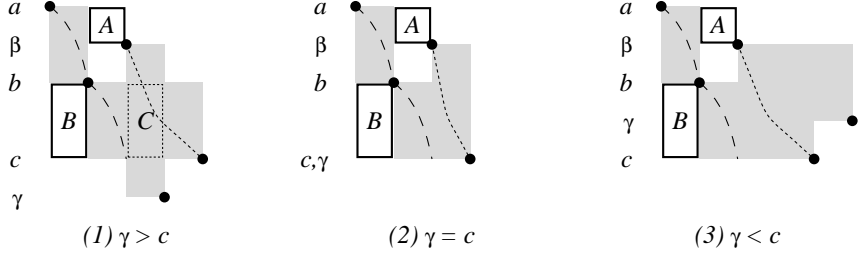
    \myfig{xtabc3}{We have displayed the case of $\beta < b$, but the
      proof of Step \ref{item:4} also holds for $\beta > b$.}
    
    Assume $\gamma > c$.  If $\min(\dxw|_C) \geq 2$, then $\tab,\tbc$
    are link incompatible for $xt_{\beta,\gamma} \leq w$.  This
    contradicts $x\in\maxsing$.  Otherwise, by the Cross
    Lemma~\ref{lem:cross}, $\tti_{\beta,c}\in\rxw$.  Then
    $\tab,\tbc$ are link incompatible for $xt_{\beta,c} \leq w$.  This
    contradicts $x\in\maxsing$.

    If $\gamma \leq c$, then $\tti_{\beta,c}\in\rxw$ and we get a
    contradiction as above.

  \item For every $\alpha,\beta,\gamma$, we have
    $\tti_{\alpha,\beta},\tti_{\beta,\gamma}\in\rxw$.\label{item:5}\\
    Suppose $\tti_{\alpha,\beta}\not\in\rxw$.  By the definition of
    $\delxw$ and the fact that $\xti = x$, we know that there exists
    $b$ such that $\tti_{\alpha,b}\in\rxw$.  Now we can apply the
    previous step to obtain a $c$ such that $\ttbc\in\rxw$.  Note
    that by Step \ref{item:6}, $t_{\alpha,b}$ and $t_{b,c}$ are link
    incompatible.  So our situation is as depicted as in
    \figref{fig:xtabc4}.

    \myfig{xtabc4}{We have displayed the case of $b > \beta$.}
    
    Using the logic of the previous step, we see that $\min(\dxw|_A) =
    0$ contradicts $x\in\maxsing$.  Hence $\tti_{\alpha,\beta}\in\rxw$
    as desired.

    The argument for showing $\tti_{\beta,\gamma}\in\rxw$ is analogous.

  \item We can shade our diagram as in \figref{fig:xtabclen1}.1.\\
    This follows immediately from the previous four steps.

  \item $\dxw \leq 1$ in \figref{fig:xtabclen1}.1.\label{item:7}\\
    We start by showing that if $a \neq \alpha$ and $b \neq \beta$,
    then $\ttab,\tti_{\alpha,\beta}\in\rxw$ implies
    $\ttalbt\in\qabxw$.

    By Steps \ref{item:6},\ref{item:4} and \ref{item:5},
    $\ttbc\in\qabxw$ for some $c$ (see \figref{fig:xtabc5}).

    \myfig{xtabc5}{} 
    Suppose $\ttalbt\not\in\qabxw$ --- i.e., $\dxw \geq 2$ on region A.
    Then $t_{a,b}$ and $\tbc$ are link incompatible for $x\talbt
    \leq w$.  This contradicts $x\in\maxsing$.
    
    A similar argument can be used to show that if $b \neq \beta$ and
    $c \neq \gamma$, then $\ttbc,\tti_{\beta,\gamma}\in\rxw$ implies
    $\tti_{\beta,\gamma}\in\qbcxw$.  

    The claim of $\dxw \leq 1$ then follows by inspection from these
    two facts.
  \item Condition 1 in \propref{prop:xtabc} holds, namely $l \geq 2$.
    
    Recall that in Step \ref{item:5} we showed that
    $\ttab,\ttbc\in\rxw$ for any $a,b,c$.  This implies that we can
    shade $\digxt$ as in \figref{fig:xtabclen1}.1.
    \myfig{xtabclen1}{}
  
    By Step \ref{item:6}, $t_{k,k+1}$ and $t_{k+1,k+l+m}$ are link
    incompatible.  This implies that $\min(\dxw|_A) = 0$.  Similarly,
    $\min(\dxw|_B) = 0$.  It then follows from our explicit
    description of $x$ that the values of $w(i)$ for $i = 1,k+1$
    are as shown in \figref{fig:xtabclen1}.2.  Arguing with $\dxw'$
    (see \ref{eq:dualdxw}) and region B, we see that $w(i)$ for
    $i=k+l,k+l+m$ is as shown in the same figure.  But this means that
    $w^{-1}(k+l) = 1$ and $w^{-1}(k+1) = k+l+m$.  This can only
    happen if $l > 1$.
  
  \item $w$ is as stated in Condition 3.
    
    Step \ref{item:7} tells us that we can conclude that $\dxw =1$ on
    all shaded areas of \figref{fig:xtabclen1}.  Therefore, $w(i) =
    k + 2 - i$ for $2 \leq i \leq k$.  A similar argument to that in
    Step \ref{item:7} shows that $w(i) = 2(k+l) + m - i$ for $k+l+1
    \leq i < k+l+m$.  So we need only investigate the values of
    $w(i)$ for $k+1 < i < k+l$.  To do this, assume that $w(i) =
    x(i)$ for $k+1 < i \leq j$ for some $j$ with $k + 1 \leq j < k
    + l - 1$.  Then, as in \figref{fig:xtabclen2}.1, we see that $\dxw
    = 0$ on region B.

    \myfig{xtabclen2}{}  
  
    Therefore, $t_{k,k+1}$ and $t_{k+1,j+2}$ are link incompatible for
    $xt_{j+2,k+l+m} \leq w$.  This contradicts $x\in\maxsing$.  Hence
    $w(i) = x(i)$ for all $i$ with $k+1 < i < k+l$.  So 
    \begin{equation}
      \begin{split}
        w = [k+l,k,\ldots,&2,k+l+m,k+l-1,\ldots,k+2,\\
                             &1,k+l+m-1,\ldots,k+l+1,k+1].      
      \end{split}
    \end{equation}
    as desired.
  
  \item Condition 2 in \propref{prop:xtabc} holds.
    
    We need to show that if $l > 2$ then $k,m = 1$.  So assume $k >
    1$.  By Steps \ref{item:5} and \ref{item:7}, $t_{1,k+1}$ and $t_{2,k+2}$
    are patch incompatible reflections for $xt_{k+l,k+l+m} \leq w$
    (see \figref{fig:xtabclen2}.2).  This contradicts $x\in\maxsing$.
    The argument showing that $m = 1$ is analogous.
  \end{enumerate}
  \renewcommand{\labelenumi}{\theenumi}

  This completes the proof of \propref{prop:xtabc}. 
\end{proof}

\section{Maximal singularity of candidates}
\label{sec:maxsing}

We now finish the proof of \thmref{thm:main} by showing that the restrictions
we have discovered for $\xti$ in Propositions \ref{prop:xtab} and
\ref{prop:xtabc} are sufficient to show that these points
correspond to \msp's in the appropriate Schubert variety.  This task
consists of two steps:
\begin{enumerate}
\item Show that the points $x$ are singular points.
\item Show that any cover of $x$ that is still below $w$ is a smooth
  point.  
\end{enumerate}

So that we can describe $\maxsing$ succinctly, we introduce the
following notation:
\begin{dfn}
  For $k,m \geq 2$, define
  \begin{align}
    x_{k,m} &= [k,\ldots,1,k+m,\ldots,k+1],\label{eq:7}\\
    w_{k,m} &= [k+m,k,\ldots,2,k+m-1,\ldots,k+1]\label{eq:8}. 
  \end{align}
  For $k,m \geq 1$ and $l \geq 2$, define 
  \begin{align}
    x_{k,l,m} &= [k,\ldots,1,k+l,\ldots,k+1,k+l+m,\ldots,k+l+1],\label{eq:9}\\
    \begin{split}
    w_{k,l,m} &= [k+l,k,\ldots,2,k+l+m,k+l-1,\ldots,k+2,\\
                &\qquad \qquad 1,k+l+m-1,\ldots,k+l+1,k+1].      
    \end{split}
  \end{align}
\end{dfn}

\begin{thm}[Rephrasing of \thmref{thm:main}]\label{thm:rephrase}
  $x$ is an \msp\ of $X_w$ if and only if
  \begin{enumerate}
  \item $t\in\rxw \text{ implies } l(xt) = l(x) + 1$.\label{eq:lencond}
  \item
    \begin{enumerate}
    \item For some $k,m \geq 2$, we have $\xti = x_{k,m}$ and $\wti = w_{k,m}$\label{item:8}\\
      or
    \item For some $k,m \geq 1$, $l = 2$ or $k = m = 1, l \geq 2$,
      we have $\xti = x_{k,l,m}$ and $\wti = w_{k,l,m}$. \label{item:9}
    \end{enumerate}
  \end{enumerate}
\end{thm}

\begin{proof}
  \propref{prop:msplo} tells us that Condition 1 is necessary.
  Propositions \ref{prop:xtab} and \ref{prop:xtabc} tell us that
  Conditions \ref{item:8} and \ref{item:9} are necessary.  So all we
  need to show is sufficiency.

  Let $t$ be a reflection such that $x < y = xt \leq w$.  As $\pht$ is
  injective, to calculate $\#\ryw$ from $\rxw$ we need only count how
  many reflections in $\rxw$ are not in the image of $\pht$.  Note by
  \propref{prop:rxtwt} that $\#\rxw = \#\rxtwt$ and $l(w) - l(x) =
  l(\wti) - l(\xti)$.

  \myfig{wtiab1}{}
  Consider first the case shown in \figref{fig:wtiab1} of two decreasing
  sequences for $\xti$.  Note that
  \begin{align}
    l(w_{k,m}) &= \binom{k}{2} + \binom{m}{2} + k + m - 1,\label{eq:lens}\\
    l(x_{k,m}) &= \binom{k}{2} + \binom{m}{2} \text{ and }\label{eq:lens2}\\
    \#\rxw &= k\cdot m.\label{eq:16}
  \end{align}
  Hence,
  \begin{equation}\label{eq:17}
    l(w) - l(x) - \#\rxw = k + m - km - 1.
  \end{equation}
  Since $k,m \geq 2$, (\ref{eq:17}) is negative.  So by
  \thmref{fact:basic}, $e_x$ is a singular point of $X_w$.  
  
  To prove that it is a maximal singular point, we consider some
  $\tab\in\rxw$ and let $y = x\tab$.
  \myfig{wtiab2}{$y = x\tab$.} 
  Then, viewing \figref{fig:wtiab2}, it is easily seen that $\#\ryw =
  (k-1) + (m-1) = k + m - 2$.  Since $l(y) = l(x) + 1$, by
  \thmref{fact:basic} and (\ref{eq:17}), $y$ is a smooth point of
  $X_w$.  Since $y$ was chosen as an arbitrary cover of $x$, $x$ is an
  \msp\ for $w$.

  \myfig{wtiabc1}{} 
  Now we prove the case shown in \figref{fig:wtiabc1} of three
  decreasing sequences for $\xti$.  Note that
  \begin{align}
    l(w_{k,l,m}) &= \binom{k}{2} + \binom{l}{2} + \binom{m}{2} + k + m
    + 2(l - 2) + 1,\label{eq:llens}\\
    l(x_{k,l,m}) &= \binom{k}{2} + \binom{l}{2} + \binom{m}{2} \text{
    and }\label{eq:llens2}\\ 
    \#\rxw &= l(k+m).
  \end{align}
  Hence,
  \begin{equation}\label{eq:14}
    l(w) - l(x) - \#\rxw = (l-1)\left( 1 + \frac{l-2}{l-1} - k -
    m\right).  
  \end{equation}
  
  Since $k,m \geq 1$ and $l \geq 2$, (\ref{eq:14}) is negative.  So by
  \thmref{fact:basic}, $x$ is a singular point of $X_w$.  To prove it
  is an \msp\ for $w$, as above we consider some $\tab\in\rxw$
  and let $y = x\tab$.  We have $l(w) - l(y) = k + m + 2(l-2)$.
  
  Viewing \figref{fig:wtiabc2}, it is clear that
  $\#\ryw = (k - 1) + (l - 1) + m(l-1)$.
  \myfig{wtiabc2}{} 
  If $l = 2$, then $\#\ryw = k + m = l(w) - l(y)$.  If $l > 2$, then
  by \ref{item:hmm} of \propref{prop:xtabc}, we have that $k = m = 1$,
  and $\#\ryw = 2(l-1) = l(w) - l(y)$.  So, in either case, $y$
  is a smooth point of $X_w$.

  So in both cases, $x$ is an \msp\ of $X_w$ as claimed.
\end{proof}

This completes the proof of \thmref{thm:rephrase}.  It is easy to
check that the above formulation is equivalent \thmref{thm:main}.
(Note, however, that the values of $k,l,m$ in the statement of
\thmref{thm:rephrase} differ from those used in \thmref{thm:main}.)


\section{Lakshmibai-Sandhya Conjecture}\label{sec:lak}

Let $w=[w(1),\ldots,w(n)]\in\fsn$.  Define $E_{w}$ to be the set of
all $x=[x(1),\ldots,x(n)]$ satisfying the following conditions:

\begin{enumerate}
\item There exist $i<j<k<l$ and $i' < j' < k' < l'$ such that (as sets)
\begin{equation*}
  \{w(i), w(j),w(k), w(l) \}=\{x(i'), x(j'), x(k'), x(l')\}.
\end{equation*}

\item One of the following holds:
\begin{enumerate}
\item $\fl_{ijkl}(w)=3412$ and $\fl_{i'j'k'l'}(x)=1324$.\label{item:2}  
\item $\fl_{ijkl}(w)=4231$ and $\fl_{i'j'k'l'}(x)=2143$.\label{item:3}
\end{enumerate}

\item Using the notation of \secref{sec:prelim}:
  \begin{enumerate}
  \item If \ref{item:2} holds, then set 
    \begin{align*}
      \underline{w} &= \ra^w_{ijkl}(1324)\text{ and }\\      
      \widehat{x} &= \ra^x_{i'j'k'l'}(3412).
    \end{align*}
  \item If \ref{item:3} holds, then set 
    \begin{align*}
      \underline{w} &= \ra^w_{ijkl}(2143)\text{ and}\\
      \widehat{x} &= \ra^x_{i'j'k'l'}(4231).
    \end{align*}
  \end{enumerate}
  Then, 
  \begin{equation}\label{eq:5}
    \underline{w} \leq x \leq \widehat{x}\leq w.
  \end{equation}
\end{enumerate}

\begin{thm}[Conjecture in \cite{Lak-San}]\label{conj:A} 
  For $w \in \fsn$, the singular locus of $X_w$ is equal to
  $\cup_{x}X_x$, where $x$ runs over the maximal elements of $E_{w}$
  in Bruhat order.
\end{thm}

\begin{proof}
  We only give the argument for singular points of the type 4231
  (i.e., those described in Case 1 of \thmref{thm:main}).  The argument
  for singular points of type 3412 and 45312 is analogous.
  
  We start by proving that $\maxsing \subseteq \cup_x X_x$.
  To do this, fix some $x\in\maxsing$ (of type 4231).  We will choose
  indices $i,j,k,l$ and $i',j',k',l'$ as described in the definition
  of $E_w$, and show that (\ref{eq:5}) is satisfied for our choice of
  indices.  So, using the notation of \thmref{thm:main}, let
  \begin{equation}\label{eq:19}
    \alpha_1 < \beta_1 < \beta_2 < \cdots < \beta_{k-1} < \alpha_2 <
    \alpha_3 < \cdots < \alpha_m < \beta_k
  \end{equation}
  correspond to a type 4231 pattern in $w$.  Then set 
  \begin{align}
    i &= \alpha_1,\quad j = \beta_1,\quad
    k = \alpha_{m-1}\ \ \ \text{ and } l = \beta_k,\\
    i' &= \alpha_1,\quad j' = \beta_{k-1},\quad
    k' = \alpha_2\ \ \ \text{ and } l' = \beta_k.
  \end{align}
  Now, recall from \lemref{lem:bdf} that $u\leq v$ if and only if
  $\duv$ is everywhere non-negative.  But then (\ref{eq:5}) follows
  from \figref{fig:ls1} along with the observation that $\duv \leq 1$
  in each of these diagrams.
  \myfig{ls1}{(1) $\underline{w},w$ (2) $x,w$ (3) $\widehat{x},w$.}
  
  Now we need to show that any $x$ satisfying (\ref{eq:5}) for some
  $w$ and some set of indices is a singular point of $X_w$.  Since $x
  \leq \widehat{x} \leq w$ by hypothesis, we see that $i,j,k,l \in
  \delxw$.  Combining this with \corref{cor:deliff} and the fact that
  $\underline{w} \leq x$, we conclude that
  \begin{equation}
    i \leq i',l' \leq l,\quad j\leq j' \leq l,\quad i \leq k' \leq k.
  \end{equation}
  By hypothesis, we also know that $i' \leq j' \leq k' \leq l'$.
  Finally, recall that
  \begin{equation}\label{eq:20}
    \{w(i), w(j),w(k), w(l) \}=\{x(i'), x(j'), x(k'), x(l')\}.
  \end{equation}
  \myfig{ls2}{}
  One possible configuration of the points of (\ref{eq:20}) is shown
  in \figref{fig:ls2}.1.  Note that (\ref{eq:5}) implies that $\dxw
  \equiv 1$ on region $A \cup B \cup C$.  Hence, $t_{i',k'},t_{j',l'}$
  are patch incompatible for $x\leq w$.  If $l(xt_{j',l'}) = l(x) +
  1$, then we can conclude that $x$ is a singular point of $X_w$ by
  \factref{fact:sing}.  Otherwise, consider \figref{fig:ls2}.2.  Pick
  the index $p$ as small as possible such that $\pt_x(p) \in D$.  Then
  $t_{i',k'}, t_{j',p}$ are patch incompatible for $xt_{p,l'} \leq
  w$.  But then $x$ is a singular point $X_w$ as $l(xt_{p,l'}t_{j',p})
  = l(xt_{p,l'}) + 1$ by construction.
  
  Carrying out the analogous arguments for the 3412 and 45312 type
  singularities (using link incompatible reflections) completes the
  proof of Conjecture \ref{conj:A}.
\end{proof}

\section{Kazhdan-Lusztig polynomials at elements of $\maxsing$}
\label{sec:klp}

The determination of $\maxsing$ has applications to the study of
multiplicities in Verma modules through the Kazhdan-Lusztig
polynomials.  These polynomials, in the type A case, are indexed by
two permutations $x,w$ (for properties of these polynomials, see
\cite{Hum}).  A result of Kazhdan and Lusztig \cite{K-L1} is that
$P_{x,w} = 1$ if and only if $e_x$ is a smooth point of $X_w \subseteq
SL(n)/B$.  Furthermore, by a conjecture of Kazhdan and Lusztig
(\cite{K-L1}), proved independently by Beilinson-Bernstein
\cite{BeilBern} and Brylinski-Kashiwara \cite{BryKash}, $P_{x,w}(1)$
gives the multiplicity of an irreducible module associated to $w$ in
the Verma module associated to $x$.

\thmref{thm:main} gives explicit conditions for this Verma module
multiplicity to be greater than 1.  In this section, we calculate the
$P_{x,w}$ for $x\in\maxsing$.  Setting $q=1$ therefore yields the
exact multiplicities in these cases.

There are some explicit descriptions of the $P_{x,w}$ in special
cases.  Lascoux and Sch\"utzenberger \cite{LS11} and Zelevinsky
\cite{Zel} combine to give a small resolution of $X_w$ and
corresponding formulas for $P_{x,w}$ when $w$ is a Grassmannian
permutation.  Lascoux \cite{Lascoux95} has extended this result to
vexillary permutations.  Brenti \cite{Brenti94,Brenti97,Brenti98} has
given several beautiful, general, alternating sum formulas.  Finally,
formulas have been calculated in several specific classes (e.g.,
\cite{Bren-Sim,hexagon,polo}).  In particular, \thmref{thm:kl45312} is
proved in \cite{ssv} and \thmref{thm:kl3412} is proved in \cite{LS11},
but both are only proved in the case where $\xti = x$ and $\wti = w$.

A result of Polo, \cite{polo}, states that every polynomial in
$\bbN[q]$ with constant term 1 can be realized as a Kazhdan-Lusztig
polynomial in $\fsn$ for some $n$.  However, as we will see in the
below three theorems, the Kazhdan-Lusztig polynomials at elements of
$\maxsing$ are of very limited forms.

For pairs of permutations $x,w\in\fsn$, we can define the
Kazhdan-Lusztig polynomials by the following properties:
\begin{enumerate}
\item $P_{x,w} = 0$ if $x\not\leq w$.
\item $P_{x,w} = 1$ if $x\leq w$ and $l(w) - l(x) \leq 2$.
\item $\deg(P_{x,w}) \leq \frac{1}{2}(l(w)-l(x)-1)$.
\item If $s \in \simref$ such that $ws < w$ then
  \begin{equation}\label{eq:11}
    P_{x,w} = q^c P_{x,ws} + q^{1-c} P_{xs,ws} - \sum_{x \leq z < ws, zs <
      z} \mu(z,ws) q^{(l(w)-l(z))/2} P_{x,z}.
  \end{equation}
  where $\mu(z,ws)$ is the coefficient of $q^{(l(ws)-l(z)-1)/2}$ in
  $P_{z,ws}$ and $c = 1$ if $xs < x$, $c = 0$ if $xs > x$.
\end{enumerate}

\begin{lem}\label{lem:tilkl}
  If $i\not\in\delxw$, then $P_{x,w} = P_{\xbi,\wbi}$
\end{lem}
\begin{proof}
  Fix $x < w$ and pick some $i\not\in\delxw$.  We know by
  \propref{prop:lblf} that $x(i) = w(i)$.  By \corref{cor:delxy}, this
  implies that if $x \leq z \leq w$ for some $z$, then $z(i) = x(i) =
  w(i)$.  
  
  With these facts, the result then follows easily by induction on
  $l(w)$ using (\ref{eq:11}).  (Note that our base case of $l(w) = 1$
  is trivial.)
\end{proof}

\begin{cor}\label{cor:pxw}
  $P_{\xti,\wti} = P_{x,w}$.
\end{cor}

As it will be used repeatedly in upcoming arguments, for reference we
state the following fact \cite[Cor. 7.14]{Hum}:
 \begin{fact}\label{fact:pxwsim}
   For $s,s'\in \simref$, $ws < w$, $s'w < w$, then $P_{x,w} =
   P_{xs,w} = P_{s'x,w}$.
 \end{fact}
 We are now ready to calculate $P_{x,w}$ for $x\in\maxsing$.  By
 \thmref{thm:rephrase} and \corref{cor:pxw}, it is enough to calculate
 $P_{x,w}$ for the pairs $x_{k,m},\ w_{k,m}$ and $x_{k,l,m},\ 
 w_{k,l,m}$.

\begin{thm}[4231-type singularities] 
  For $k,m \geq 2$, 
  \begin{equation}
    P_{x_{k,m},w_{k,m}} = 1 + q + \cdots + q^{\min(k-1,m-1)}.
  \end{equation}
\end{thm}
\begin{proof}
  We apply induction on $k + m$.  The case of $k = m = 2$ can be
  checked from (\ref{eq:11}).  We assume then, without loss of
  generality, that $k \geq 3$.  Also, for brevity we will often write
  $x$ and $w$ in place of $x_{k,m}$ and $w_{k,m}$, respectively.
  
  \myfig{kl4231}{}
  Consider $s = s_1$.  The pairs $x,w$ and $x,ws$ and $xs,ws$ are shown
  in \figref{fig:kl4231}.  We claim that 
  \begin{equation}\label{eq:18}
  P_{xs,ws} = P_{x,ws} = 1 + q + \cdots + q^{\min(k-2,m-1)}.
  \end{equation}
  The second equality follows from \figref{fig:kl4231}.2,
  \corref{cor:deliff}, \lemref{lem:tilkl} and the induction
  hypothesis.  To obtain the first equality, we notice that $s_{k-1}xs
  = x$.  Then, since $s_{k-1}ws < ws$, (\ref{eq:18}) follows from
  \factref{fact:pxwsim}.

  Substituting this information into (\ref{eq:11}), we obtain
  \begin{equation}\label{eq:13}
    P_{x,w} = (1+q)\left(1 + q + \cdots + q^{\min(k-2,m-1)}\right) - 
    \nts \nts \sumsb{x \leq z < ws \\ zs < z}\nts\nts \mu(z,ws) 
    q^{\frac{l(w)-l(z)}{2}} P_{x,z}.
  \end{equation}
  
  Now we need to investigate the possible terms in the sum of
  (\ref{eq:13}).  We will first determine which $z$ with $l(z) < l(ws)
  - 1$ can contribute to (\ref{eq:13}).  By \corref{cor:delxy},
  \corref{cor:deliff} and \lemref{lem:tilkl}, if $x \leq z < ws$ then
  $P_{z,ws} = P_{\zbio,\wsbio} = P_{x_{k-1,m},w_{k-1,m}}$.  To have
  $\mu(z,ws) > 0$ while $l(z) < l(ws) - 1$, we need $\deg(P_{z,ws}) >
  0$.  Since $x_{k-1,m}$ is the only\ \msp\ for $w_{k-1,m}$, this
  tells us that $\zbio \leq x_{k-1,m}$.  But by \factref{fact:pxwsim},
  $P_{e,w_{k-1,m}} = P_{x_{k-1,m},w_{k-1,m}}$, so this implies that we
  must have $\zbio = x_{k-1,m}$ by the degree bound.  It is then easy
  to check from (\ref{eq:lens}) and (\ref{eq:lens2}) that $\mu(x,ws) >
  0$ if and only if $k-1 = m$.

  We now split into two cases depending on the relative values of
  $k$ and $m$.
  \begin{enumerate}
  \item $k-1 \neq m$.\\
    Since $k-1 \neq m$, we know by the previous paragraph that the
    only $z$ that can contribute to the sum in (\ref{eq:13}) are those
    with $l(z) = l(ws) - 1$.  Furthermore, since we are summing only over $z$
    for which $zs < z$, the only possibility is $z = wss_2$.  In this
    case, by induction, along with Corollaries \ref{cor:deliff} and
    \ref{cor:pxw}, $P_{x,wss_2} = 1 + q + \cdots + q^{\min(k-3,m-1)}$.
    So
    \begin{equation}
    \begin{split}
      P_{x,w} &= (1+q)(1 + q + \cdots + q^{\min(k-2,m-1}) - 
      q(1+ q \cdots + q^{\min(k-3,m-1)})\\ 
      &= 1 + q + \cdots + q^{\min(k-1,m-1)}.
    \end{split}
    \end{equation}
    
  \item $k-1 = m$.\\
    $z = wss_2$ will contribute as in the previous case.  However,
    from our discussion above, $z = x_{k-1,m}$ will also contribute.
    For $z = x_{k-1,m}$, we have $P_{x,z} = 1$ and $q^{(l(w)-l(z))/2}
    = q^{k-1}$.  Plugging this term into (\ref{eq:13}), along with the
    term coming from $z = wss_2$, we get
    \begin{equation}
      \begin{split}
        P_{x,w} &= (1+q)(1 + q + \cdots + q^{k-2}) -  q(1+ q
        \cdots + q^{k-3}) - q^{k-1}\\
        &= 1 + q + \cdots + q^{k-2} = 1 + q + \cdots + q^{\min(k-1,m-1)}.
      \end{split}
    \end{equation}
  \end{enumerate}
\end{proof}

\begin{thm}[3412-type singularities]\label{thm:kl3412}  
  For $k,m \geq 1$, 
  \begin{equation}
    P_{x_{k,2,m},w_{k,2,m}} = 1 + q.
  \end{equation}
\end{thm}
\begin{proof}
  We apply induction on $k + m$.  The case of $k = m = 1$ is easy to
  check using (\ref{eq:11}).  For brevity, we abbreviate $x_{k,2,m}$
  and $w_{k,2,m}$ by $x$ and $w$, respectively.
  
  \myfig{kl3412}{}
  
  Let $s = s_{k+1}$.  Now, as seen in \figref{fig:kl3412}.2, $ws$ is
  3412- and 4231-avoiding, hence smooth (\cite{Lak-San}).  Therefore,
  $P_{y,ws} = 1$ for all $y\leq ws$.  Clearly $x,xs \leq ws$.  Thus
  the first two terms of (\ref{eq:11}) together contribute $1 + q$.
  We now show that the sum in (\ref{eq:11}) is empty.
  
  Since $P_{y,ws} = 1$ for all $y\leq ws$, $\mu(z,ws) > 0$ implies
  that $l(z) = l(ws) - 1$.  But, as seen from \figref{fig:kl3412}.2,
  no $z$ satisfying this length condition can satisfy the additional
  constraint of $zs_{k+1} < z$.  This proves the theorem.
\end{proof}

\begin{thm}[45312-type singularities]\label{thm:kl45312} 
  For $l \geq 2$, 
  \begin{equation}
    P_{x_{1,l,1},w_{1,l,1}} = 1 + q^{l-1}.
  \end{equation}
\end{thm}
\begin{proof}
  We apply induction on $l$.  The case of $l = 2$ is covered by
  \thmref{thm:kl3412}.  So we assume $l\geq 3$.  For clarity, we
  abbreviate $x_{1,l,1}$ and $w_{1,l,1}$ by $x$ and $w$, respectively.
  
  \myfig{kl45312}{}
  
  In \figref{fig:kl45312}, we depict the pairs $x,\ w$ and $x,\ ws_2$
  and $xs_2,\ ws_2$.  We claim that the first two terms in
  (\ref{eq:11}) contribute $(1+q)(1+q^{l-2})$.  First consider the
  pair $x,\ ws_2$.  Since $ws_2s_1 < ws_2$, by the induction
  hypothesis, \corref{cor:deliff} and \lemref{lem:tilkl}, we see that
  $P_{x,ws_2} = P_{xs_1,ws_2} = 1 + q^{l-2}$.  Now consider the pair
  $xs_2,\ ws_2$.  Since $s_lws_2 < ws_2$ and $ws_2s_1 < ws_2$, it
  follows that $P_{xs_2,ws_2} = P_{s_lxs_2s_1,ws_2}$.  But since
  $s_lxs_2s_1 = xs_1$, we get that $P_{xs_2,ws_2} = 1 + q^{l-2}$ also.
  Plugging this information into (\ref{eq:11}), we can write
  \begin{equation}\label{eq:12}
    P_{x,w} = 1 + q^{l-2} + q + q^{l-1} - \sumsb{x \leq z < ws\\ zs_2 <
      z} \mu(z,ws_2) q^{\frac{l(w)-l(z)}{2}} P_{x,z}.
  \end{equation}
  
  Now we check which $z$ will appear in the sum in (\ref{eq:12}).
  First note that $xs_1$ is the unique \msp\ for $ws_2$.  By
  induction, $P_{xs_1,ws_2} = 1 + q^{l-2}$.  By \factref{fact:pxwsim},
  $P_{e,ws_2} = P_{xs_1,ws_2}$.  Hence, the only $z$ such that $l(z) <
  l(ws_2) - 1$ and $\deg(P_{z,ws_2})$ is maximized is $z = xs_1$.
  However, $xs_1s_2 > xs_1$, so $xs_1$ does not appear in the sum.  So
  the only possible terms in the sum are those with $l(z) = l(ws_2) -
  1$.  From \figref{fig:kl45312}.2, we see that $z = ws_2s_3$ is the
  only $z$ satisfying both this length condition and $zs_2 < z$.
  Using \factref{fact:pxwsim}, \lemref{lem:tilkl} and the induction
  hypothesis, one can check that $P_{x,ws_2s_3} = 1 + q^{l-3}$.
  Hence, the sum in (\ref{eq:12}) contributes $-q-q^{l-2}$.
  Simplifying, we see that $P_{x,w} = 1 + q^{l-1}$ as claimed.
\end{proof}

\section{Examples calculating $\maxsing$}
\label{sec:ex}

\begin{exmp}
  Using \thmref{thm:main}, in \figref{fig:example} we compute the
  singular locus 
  \begin{equation}
  \maxsing = X_{[48376512]} \cup X_{[64387512]} \cup
  X_{[46587312]} \cup X_{[68174325]}
\end{equation}
of $X_w$ where $w = [6,8,4,7,5,3,1,2]$.
  \myfig{example}{}
\end{exmp}

\begin{remark}
  The cardinality of the set $\maxsing$ may be $O(n^4)$.  This is the
  case, for example, when 
  \begin{equation}
    w = [k+1,\ldots,k+l,1,\ldots,k].    
  \end{equation}
  Then $\#\maxsing = \binom{n/2}{2}^2$.
\end{remark}

\begin{exmp}
  Using a computer it is easy to calculate, for example, that for 
  \begin{equation}
    w = [17,6,2,15,12,11,3,8,16,7,14,5,13,9,10,1,4],
  \end{equation}
  $\#\maxsing = 29$.
\end{exmp}

\section{Patterns indexing $\maxsing$}\label{sec:ind}

Which 4231 or 3412 patterns lead to elements in $\maxsing$?  We can
describe these patterns by taking all 4231 and 3412 patterns in $w$
and removing certain ``useless patterns'' contained in larger patterns
of length 5 or 6.  For example, if $w=[52341]$, the pattern 5241 will
be useless since the shaded region it defines is not empty.  We
describe the useless patterns in the following way.  For each pattern
of length 5 or 6 in the left hand column of \eqref{table:useless},
remove the corresponding pattern in the right hand column.

\begin{equation}\label{table:useless}
\begin{matrix}
(5 2 3 4 1) & (5 2 4 1) \\
(5 2 4 3 1) & (5 2 4 1) \\
(5 3 2 4 1) & (5 2 4 1) \\
(5 3 4 2 1) & (5 3 4 1) \\
(5 4 2 3 1) & (5 2 3 1) \\
(3 5 4 1 2) & (3 5 1 2) \\
(4 3 5 1 2) & (4 5 1 2) \\
(4 5 1 3 2) & (4 5 1 2) \\
(4 5 2 1 3) & (4 5 1 3) \\
(6 3 5 2 4 1) & (6 3 4 1) \\
(5 6 3 4 1 2) & (5 6 1 2) \\
(5 2 6 4 1 3) & (5 6 1 3) \\
(4 6 3 1 5 2) & (4 6 1 2) 
\end{matrix}
\end{equation}

The remaining ``useful patterns'' all index a unique component in
$\maxsing$.  For example, if $w = [7432651]$ then $\maxsing$ has only
one element namely $x=[4321765]$ and this element would be indexed by
$7251$.  This example corresponds to the shape in
Figure~\ref{fig:mainthm1}.

It would be interesting to know the distribution of the various sizes
of $\maxsing$ for all $w \in \fsn$ for large $n$.

\bibliography{gen}
\end{document}